\magnification=\magstep1
\input amstex
\documentstyle{amsppt}

\define\defeq{\overset{\text{def}}\to=}
\define\ab{\operatorname{ab}}
\define\pr{\operatorname{pr}}
\define\Gal{\operatorname{Gal}}
\define\diag{\operatorname{diag}}
\define\Hom{\operatorname{Hom}}
\define\sep{\operatorname{sep}}
\def \isom {\overset \sim \to \rightarrow}
\define\cor{\operatorname{cor}}
\define\Pic{\operatorname{Pic}}
\define\Br{\operatorname{Br}}

\define\tame{\operatorname{tame}}

\define\id{\operatorname{id}}
\define\cn{\operatorname{cn}}
\define\inn{\operatorname{inn}}
\define\Ker{\operatorname{Ker}}

\def \c{\operatorname {c}}
\def \Aut{\operatorname {Aut}}
\def \Out{\operatorname {Out}}

\def \out{\operatorname {out}}

\def \Sec{\operatorname {Sec}}
\def \Inn{\operatorname {Inn}}
\def \Sp{\operatorname {Sp}}

\def \char{\operatorname{char}}
\def \Br{\operatorname{Br}}

\define\Primes{\frak{Primes}}
\NoRunningHeads
\NoBlackBoxes
\topmatter

\title
Around the  Grothendieck anabelian section conjecture
\endtitle

\author
Mohamed Sa\"\i di
\endauthor

\abstract 
This paper is around the topics I discussed in the lecture I gave at the Isaac Newton Institute in Cambridge, July 2009, in the Introductory Workshop.
This paper can be read as a companion to my paper [Sa\"\i di], where detailed proofs can be found.
\endabstract
\toc

\subhead
\S 0. Introduction 
\endsubhead

\subhead
\S 1. Generalities on Arithmetic Fundamental Groups and Sections
\endsubhead

\subhead
\S 2. Grothendieck Anabelian Section Conjecture
\endsubhead

\subhead
\S 3. Good Sections of Arithmetic Fundamental Groups
\endsubhead

\subhead
\S 4. The Cuspidalisation of Sections of Arithmetic Fundamental groups
\endsubhead

\subhead
\S 5. Applications to the Grothendieck Anabelian Section Conjecture
\endsubhead

\subhead
\S 6. On a Weak Form of the $p$-adic Grothendieck Anabelian Section Conjecture
\endsubhead

\endtoc
\endtopmatter

\document

\subhead
\S 0. Introduction
\endsubhead
This, mostly expository, paper is around the topic of the Grothendieck anabelian section conjecture.

This conjecture predicts that splittings, or sections, of the exact
sequence of the arithmetic fundamental group 
$$1\to \pi_1(\overline X)\to \pi_1(X)\to G_k\to 1$$
of a proper, smooth, and hyperbolic curve $X$, all arise from decomposition subgroups associated to rational points of $X$, over
a base field $k$ which is finitely generated over the prime field $\Bbb Q$ (cf. $\S2$).

A birational version of this conjecture predicts that splittings of the exact sequence of absolute Galois groups
$$1\to \Gal (K_X^{\sep}/\bar k.K_X)\to \Gal (K_X^{\sep}/K_X)\to G_k\to 1,$$
where $K_X$ is the function field of the hyperbolic curve $X$, all arise from rational points of $X$ in some precise way (cf. $\S2$), 
under the above assumption on the base field $k$ (cf. $\S2$).

This conjecture is one of the main topics in anabelian geometry. 
It establishes a dictionary between profinite group theory, arithmetic geometry, and Diophantine geometry.

This conjecture is still widely open. Only some examples are treated in the literature.

The most complete achievement around this conjecture, is the proof by Koenigsmann, and Pop, that the  birational version 
of this conjecture holds true over $p$-adic local fields (cf. $\S5$).

The main issue in investigating the section conjecture is the following. How can one produce a rational point $x\in X(k)$, from a splitting
of the exact sequence of $\pi_1(X)$?

This issue is completely settled in the case of birational sections over $p$-adic local fields by Koenigsmann, and Pop.
In this case, one produces a rational point by resorting to a local-global principle for Brauer groups of fields of transcendence 
degree $1$ over $p$-adic local fields, that was proven by Lichtenbaum, in the finitely generated case, and by Pop in general. 

It is not clear for the time being, at least to the author, how to settle the above issue in the case of sections of $\pi_1(X)$.

The main idea we would like to advocate in this paper is to reduce the solution of the section conjecture to the solution of its birational version.

We introduce, and investigate, the theory of cuspidalisation of sections of arithmetic fundamental groups
for this purpose (cf. $\S4$). 

The main aim of this theory is, starting from a section 
$$s:G_k\to \pi_1(X)$$ 
of the exact sequence of $\pi_1(X)$, to 
construct a section 
$$\tilde s:G_k\to \Gal (K_X^{\sep}/K_X)$$ 
of the exact sequence of $\Gal (K_X^{\sep}/K_X)$, which lifts the section $s$, i.e.
which inserts into a commutative diagram:
$$
\CD
G_k @>{\tilde s}>> \Gal (K_X^{\sep}/K_X)\\
@V{\id}VV   @VVV \\
G_k @>{s}>>  \pi_1 (X)\\
\endCD
$$
where the right vertical homomorphism is the natural one. 

We have a natural exact sequence
$$1\to I_X\to \Gal (K_X^{\sep}/K_X)\to \pi_1(X)\to 1,$$
where $I_X$ is the normal subgroup of $\Gal (K_X^{\sep}/K_X)$ generated by the inertia subgroups
at all geometric points of $X$.

We exhibit the notion of (uniformly) good sections of arithmetic fundamental groups (cf. $\S3$).
Sections which arise from rational points are good sections.

Our main result on the theory of cuspidalisation is that good sections $s:G_k\to \pi_1(X)$
of $\pi_1(X)$ can be lifted to sections 
$${\tilde s}^{c-ab}:G_k\to \Gal (K_X^{\sep}/K_X)^{\c-\ab}$$
of the natural projection $\Gal (K_X^{\sep}/K_X)^{\c-\ab}\twoheadrightarrow G_k$, where
$\Gal (K_X^{\sep}/K_X)^{\c-\ab}$ is defined by the following push out diagram:

$$
\CD
1   @>>>   I_X  @>>>   \Gal (K_X^{\sep}/K_X)  @>>> \pi_1(X)  @>>> 1  \\
@.   @VVV         @VVV         @V{\id}VV    \\
1  @>>>   I_X^{\ab}   @>>>  \Gal (K_X^{\sep}/K_X) ^{\c-\ab}  @>>>   \pi_1 (X)  @>>> 1 \\
\endCD
$$
and $I_X^{\ab}$ is the maximal abelian quotient of $I_X$, under quite general assumptions
on the field $k$, which are satisfied by number fields, and $p$-adic local fields (cf. Corollary 4.7).

As an application, we prove a (pro-$p$) version of the section conjecture over $p$-adic local fields,
for good sections of arithmetic fundamental groups, under some additional assumptions (cf. Theorem 5.4).

We state a (an unconditional) result on a semi-birational version of the $p$-adic section conjecture (cf. Theorem 5.6).

We also prove that the existence of good sections of $\pi_1(X)$, over number fields, implies the existence of degree $1$ divisor on $X$, under a finiteness assumption of the 
Tate-Shafarevich group of the jacobian of $X$ (cf. Theorem 5.8).

Finally, in $\S6$, we discuss a weak version of the section conjecture over $p$-adic local fields, which is related to the absolute anabelian geometry
of hyperbolic curves over $p$-adic local fields.

\subhead
\S 1. Generalities on Arithmetic Fundamental Groups and Sections
\endsubhead
In this section we introduce the set-up of the Grothendieck anabelian section conjecture.

\subhead {1.1}
\endsubhead
Let $k$ be a field of characteristic $p\ge 0$, and $X$ a proper, smooth, geometrically connected,
and hyperbolic algebraic curve over $k$. Let $K\defeq K_X$ be the function field of $X$, and 
$\eta $ a geometric point of $X$ above 
the generic point of $X$. Then $\eta$ determines naturally an algebraic closure $\bar k$
of $k$, a separable closure $K_X^{\sep}$ of $K_X$,
and a geometric point $\bar \eta$ of $\overline X\defeq X\times _k \bar k$. 

There exists a canonical exact sequence of profinite groups

$$1\to \pi_1(\overline X,\bar \eta)\to \pi_1(X, \eta) @>{\pr_X}>> G_k\to 1.\tag {$1.1$}$$

Here, $\pi_1(X, \eta)$ denotes the arithmetic \'etale fundamental group of $X$ with base
point $\eta$, $\pi_1(\overline X,\bar \eta)$ the \'etale fundamental group of $\overline X$ with base
point $\bar \eta$, and $G_k\defeq \Gal (\bar k/k)$ the absolute Galois group of $k$.

We will consider the following variant of the above exact sequence (1.1). 
Let
$$\Sigma \subseteq \Primes$$
be a non-empty subset of the set $\Primes$ of all prime integers. 
In the case where $\char(k)=p>0$, we will assume that $p\notin \Sigma$.

Write
$$\Delta_X\defeq \pi_1(\overline X,\bar \eta)^{\Sigma}$$
for the maximal pro-$\Sigma$ quotient of $\pi_1(\overline X,\bar \eta)$,
and
$$\Pi_X\defeq  \pi_1(X, \eta)/ \Ker  (\pi_1(\overline X,\bar \eta)\twoheadrightarrow
\pi_1(\overline X,\bar \eta)^{\Sigma})$$
for the quotient of  $\pi_1(X, \eta)$ by the kernel of the natural surjective homomorphism
$\pi_1(\overline X,\bar \eta)\twoheadrightarrow \pi_1(\overline X,\bar \eta)^{\Sigma}$, which is a normal subgroup
of $\pi _1(X,\eta)$. 

Thus, we have an exact sequence of profinite groups
$$1\to \Delta_X\to \Pi_X @>{\pr_{X,\Sigma}}>> G_k\to 1.\tag {$1.2$}$$
We shall refer to $\pi_1(X, \eta)^{(\Sigma)}\defeq \Pi_X$
as the geometrically pro-$\Sigma$ quotient of $\pi_1(X, \eta)$, or the geometrically pro-$\Sigma$
arithmetic fundamental group of $X$.

The exact sequence $(1.2)$ induces a natural homomorphism
$$\rho _{X,\Sigma}: G_k\to \Out (\Delta _X),$$
where $\Out (\Delta _X)\defeq \Aut (\Delta_X)/\Inn (\Delta_X)$ is the group of outer automorphisms of $\Delta_X$.
For $g\in G_k$, its image $\rho _{X,\Sigma}(g)$ in $\Out(\Delta _X)$
is the class of the automorphism of $\Delta_X$ obtained by
lifting $g$ to an element $\tilde g\in \Pi_X$, and let $\tilde g$ act on $\Delta _X$ by inner conjugation.

A profinite group $G$ is slim if every open subgroup of $G$ is centre-free (cf. [Mochizuki], \S 0).
If $G$ is a slim profinite group, we have a natural exact sequence
$$1\to G\to \Aut G \to \Out G \to 1,$$
where the homomorphism $G\to \Aut G$ sends an element $g\in G$ to the corresponding inner
automorphism $h\mapsto ghg^{-1}$. 

Moreover, if the profinite group $G$ is finitely generated, then the groups 
$\Aut(G)$, and $\Out(G)$, are naturally endowed with a profinite topology, and the above sequence is an exact 
sequence of profinite groups.

The following is an important property of the profinite groups $\Delta_X$, and $\Pi_X$.

\proclaim {Lemma 1.2} The profinite group $\Delta_ X$ is slim. In particular, the exact sequence
(1.2) is obtained from the following exact sequence
$$1\to \Delta _X\to \Aut (\Delta _X) \to \Out (\Delta _X)\to 1, \tag  {$1.3$}$$
by pull back via the natural continuous homomorphism 
$$\rho _{X,\Sigma}:G_k\to \Out (\Delta _X).$$ 

More precisely, we have a commutative diagram:
$$
\CD
1 @>>>  \Delta _X     @>>> \Aut (\Delta _X)     @>>> \Out (\Delta _X)  @>>> 1\\
  @.        @A{\id}AA           @AAA              @A{\rho _{X,\Sigma}}AA \\
1 @>>>  \Delta _X       @>>> \Pi_X       @>>> G_k   @>>> 1
\endCD
\tag {$1.4$}
$$
where the horizontal arrows are exact, and the right square is cartesian.

\endproclaim

\demo {Proof} Well known, see for example [Tamagawa], Proposition 1.11.
\qed
\enddemo

\subhead {1.3}
\endsubhead
Similarly, if we write
$X\times X\defeq X\times _kX,$ and $\iota : X\to X\times X$
for the natural diagonal embedding, then  the geometric point $\eta$ determines naturally (via $\iota$)
a geometric point,
which we will also denote $\eta$, of $X\times X$. 

There exists a natural exact sequence of
profinite groups
$$1\to \pi_1(\overline {X\times X},\bar \eta)\to \pi_1(X\times X, \eta) @>{\pr}>> G_k\to 1.
\tag {$1.5$}$$

Here, $\pi_1(X\times X, \eta)$ denotes the arithmetic \'etale fundamental group of $X\times X$
with base point $\eta$, which is naturally identified with the fibre product
$\pi_1(X, \eta)\times _{G_k}\pi_1(X, \eta)$, and $\pi_1(\overline {X\times X},\bar \eta)$ is the
\'etale fundamental group of $\overline {X\times X}$
with base point $\bar \eta$ ($\bar \eta$ is naturally induced by $\eta$), which is naturally
identified with the product $\pi_1(\overline X, \bar \eta)\times \pi_1(\overline X, \bar \eta)$.

Similarly, as in 1.1, we consider the maximal pro-$\Sigma$ quotient
$$\Delta_{X\times X}\defeq \pi_1(\overline {X\times X},\bar \eta)^{\Sigma}$$ 
of $\pi_1(\overline {X\times X},\bar \eta)$,
which is naturally identified with $\Delta_X\times \Delta _X$,
and the geometrically pro-$\Sigma$ quotient
$$ \Pi_{X\times X}\defeq \pi_1(X\times X, \eta)^{(\Sigma)}\defeq \pi_1(X\times X, \eta)/\Ker (\pi_1(\overline {X\times X},\bar \eta)
\twoheadrightarrow \pi_1(\overline {X\times X},\bar \eta)^{\Sigma})$$ 
of $\pi_1(X\times X, \eta)$, which is
naturally identified with $\Pi_X\times _{G_k}\Pi_X$.

Thus, we have a natural exact sequence
$$1\to \Delta_{X\times X} \to  \Pi_{X\times X} \to G_k\to 1.$$

\subhead {1.4}
\endsubhead
Our main objects of interest are group-theoretic splittings, or sections,  of the exact sequence (1.2). 

Let
$$s:G_k\to \Pi_X$$
be a continuous group-theoretic section of the natural projection 
$$\pr\defeq \pr_{X,\Sigma}:\Pi_X
\twoheadrightarrow G_k,$$
meaning that $\pr \circ s :G_k\to G_k$ is the identity homomorphism. 

We will refer to $s:G_k\to \Pi_X$, as above, as a section of the arithmetic fundamental group $\Pi_X$.

Every inner automorphism $\inn ^{g}:\Pi_X\to \Pi_X$ of $\Pi_X$ by an element $g \in \Delta _X$, gives rise to a
conjugate section $\inn ^{g}\circ s:G_k\to \Pi_X.$
We will refer to the set
$$\Cal C[s]\defeq \{\inn ^{g}\circ s:G_k\to \Pi_X\}_{g \in \Delta _X}$$
as the set of conjugacy classes of the section $s$.

\subhead
\S 2. Grothendieck Anabelian Section Conjecture
\endsubhead
We follow the same notations as in $\S1$. 

\subhead {2.1}
\endsubhead
Sections of arithmetic fundamental groups arise naturally from rational points. 

More precisely,
Let $x\in X(k)$ be a rational point. Then $x$ determines a decomposition subgroup
$$D_x\subset \Pi_X,$$
which is defined only up to conjugation by the elements of $\Delta _X$, and which maps isomorphically
to $G_k$ via the projection $\pr:\Pi_X\twoheadrightarrow G_k$. 

Hence, the subgroup $D_x\subset \Pi_X$, determines a group-theoretic
section
$$s_x:G_k\to \Pi_X$$
of the natural projection $\pr:\Pi_X \twoheadrightarrow G_k$, which is defined only up to conjugation by
the elements of $\Delta _X$.

Let $\overline {\Sec} _{\Pi_X}$ be the set of conjugacy classes of all continuous group-theoretic sections $G_k\to \Pi_X$,
of the natural projection $\pr:\Pi_X \twoheadrightarrow G_k$, modulo inner conjugation by the elements of $\Delta _X$ (cf. 1.4).

We have a natural set-theoretic map
$$\varphi_{X,\Sigma}: X(k)\to  \overline {\Sec} _{\Pi_X},$$
$$\ \ \ \ \ \ x\mapsto \varphi_{X,\Sigma} (x)\defeq [s_x],$$
where $[s_x]$ denotes the class of a section $s_x:G_k\to \Pi_X$, associated to the rational point $x$,  in $\overline {\Sec} _{\Pi_X}$.

\definition {Definition 2.2}\ Let $s:G_k\to \Pi_X$ be a continuous group-theoretic section of the natural projection
$\pr:\Pi_X \twoheadrightarrow G_k$. We say that the section $s$ is point-theoretic, or geometric, if the class $[s]$ of $s$ in
$\overline {\Sec} _{\Pi_X}$ belongs to the image of the map $\varphi_{X,\Sigma}$.
\enddefinition

The following is the main Grothendieck anabelian conjecture regarding sections of arithmetic fundamental groups. This conjecture establishes
a dictionary between purely group-theoretic sections of arithmetic fundamental groups, and rational points of hyperbolic curves 
over finitely generated fields, in characteristic $0$.

\definition {Grothendieck Anabelian Section Conjecture (GASC)} (cf. [Grothendieck])
Assume that $k$ is finitely generated over the prime field
$\Bbb Q$, and $\Sigma =\Primes$. Then the map $\varphi_X\defeq \varphi_{X,\Sigma}:
X(k)\to  \overline {\Sec} _{\Pi_X}$ is bijective. In particular, every group-theoretic section of $\Pi_X$ is point-theoretic in this case.
\enddefinition

The injectivity of the map $\varphi_X$ is well-known (cf. for example [Mochizuki1], Theorem 19.1).

So the statement of the GASC is equivalent to the
surjectivity of the set-theoretic map $\varphi_X$, i.e. that every group-theoretic section of $\Pi_X$ is point-theoretic, under the above assumptions. 

Note that a similar conjecture can be formulated over any field, and for
any non-empty set $\Sigma$ of prime integers,
but one can not expect its validity in general. 

For example, the analog of this conjecture doesn't hold over
finite fields (even if $\Sigma=\Primes$).  

Indeed, over a finite field the natural
projection $\Pi_X \twoheadrightarrow G_k$ admits group-theoretic sections, or splittings, since the profinite group
$G_k$ is free in this case. On the other hand, there are proper,
smooth, geometrically connected, and hyperbolic curves over finite fields with no rational points.

\subhead {2.3}
\endsubhead
Assume that $k$ is a number field, i.e. $k$ is a finite extension of the
field of rational numbers $\Bbb Q$. Let $v$ be a place of $k$, and denote by $k_v$ the completion of $k$ at $v$.
Write 
$$X_v\defeq X\times _k k_v.$$

Let $D_v\subset G_k$ be a decomposition group at $v$ ($D_v$ is only defined up to conjugation), which is naturally
isomorphic to the absolute
Galois group $G_{k_v}$ of $k_v$. 

By pulling back the exact sequence
$$1\to \Delta _X\to \Pi_X\to G_k\to 1,$$
by the natural injective homomorphism $D_v\hookrightarrow  G_k$, we obtain the exact sequence
$$1\to \Delta _{X_v}\to \Pi_{X_v}\to G_{k_v}\to 1.$$

Note that there exists an isomorphism $\Delta _X\isom \Delta _{X_v}$.

In particular, a group-theoretic section $s:G_k\to \Pi_X$ of the arithmetic fundamental group $\Pi_X$, induces naturally a group-theoretic 
section 
$$s_v: G_{k_v}\to \Pi_{X_v}$$ 
of $\Pi_{X_v}$, for each place $v$ of $k$. 

Moreover, if the section $s$ is point-theoretic, then the section $s_v$ is point-theoretic,
for each place $v$ of $k$, as is easily seen.

It seems quite natural to formulate an analog of the GASC over $p$-adic local fields.

\definition {A $p$-adic Version of Grothendieck Anabelian Section Conjecture (p-adic GASC)} Assume that $k$ is a 
$p$-adic local field, i.e. $k$ is a finite extension of the field $\Bbb Q_p$, for some prime integer $p$, 
and $\Sigma=\Primes$. Then the map
$\varphi_X\defeq \varphi_X\defeq \varphi_{X,\Sigma}: X(k)\to  \overline {\Sec} _{\Pi_X}$ is bijective.
\enddefinition

The map $\varphi_X$ is known to be injective in the case where $k$ is a $p$-adic local field, and
$p\in \Sigma$ (cf. [Mochizuki1], Theorem 19.1).
Thus, the statement of the $p$-adic GASC is equivalent to the surjectivity of the map $\varphi_X$ in this case.

\subhead {2.4}
\endsubhead
Next, we recall the definition of a system of neighbourhoods of a group-theoretic
section of the arithmetic fundamental group $\Pi _X$.

The profinite group $\Delta _X$ being topologically finitely generated, there
exists a sequence of characteristic open subgroups
$$...\subseteq \Delta _X[i+1]\subseteq \Delta _X[i]\subseteq...\subseteq \Delta _X[1]\defeq \Delta _X$$
(where $i$ ranges over all positive integers) of $\Delta _X$, such that 
$$\bigcap _{i\ge 1}\Delta _X[i]=\{1\}.$$

In particular, given a group-theoretic section $s:G_k\to \Pi_X$ of $\Pi_X$, we obtain open subgroups
$$\Pi _X[i,s]\defeq s(G_k).\Delta _X[i] \subseteq \Pi_X$$
(where $s(G_k)$ denotes the image of $G_k$ in $\Pi_X$ via the section $s$) of $\Pi_X$,
whose intersection coincide with $s(G_k)$, and
which correspond to a tower of finite \'etale (not necessarily Galois) covers
$$...\to X_{i+1}[s]\to X_{i}[s]\to ...\to X_1[s]\defeq X$$
defined over $k$. 

We will refer to the set $\{X_i[s]\}_{i\ge 1}$
as a system of neighbourhoods of the section $s$.

Note that for each positive integer $i$, the open subgroup $\Pi _X[i,s]$ of $\Pi_X$
is naturally identified with the geometrically pro-$\Sigma$ arithmetic \'etale fundamental group
$\pi_1(X_i[s], \eta_i)^{(\Sigma)}$
(the geometric point $\eta_i$ of $X_i[s]$ is naturally induced by the geometric point $\eta$ of $X$), and
sits naturally in the following exact sequence
$$1\to \Delta _X[i]\to \Pi_X[i,s] \to G_k\to 1,$$
which inserts in the following commutative diagram:
$$
\CD
1 @>>>  \Delta _X[i]     @>>> \Pi_X[i,s]     @>>> G_k @>>> 1\\
  @.        @VVV           @VVV              @V{\id}VV \\
1 @>>>  \Delta _X       @>>> \Pi_X       @>>> G_k   @>>> 1
\endCD
$$
where the two left vertical homomorphisms are the natural inclusions.

In particular, by the very definition of $\Pi_X[i,s]$,
the section $s$ restricts naturally to a group-theoretic section
$$s_i:G_k\to \Pi_X[i,s]$$
of the natural projection $\Pi_X[i,s]\twoheadrightarrow G_k$,
which fits into the following commutative diagram:
$$
\CD
G_k @>{s_i}>> \Pi_X[i,s] \\
@V{\id}VV     @VVV  \\
G_k @>{s}>> \Pi_X
\endCD
$$
where the right vertical homomorphism is the natural inclusion.

Thus, a section $s:G_k\to \Pi_X$ gives rise to a system of neighbourhoods $\{X_i[s]\}_{i\ge 1}$, 
and the corresponding open subgroups $\{\Pi _X[i,s]\}_{i\ge 1}$ inherit naturally, from the section $s$, 
sections $s_i:G_k\to \Pi_X[i,s]$. 

In investigating the point-theorecity of the section $s$ (cf. Definition 2.2), it is important to observe not only the section $s$, 
but rather the family of sections $\{s_i\}_{i\ge 1}$. 

In fact, a number of important properties of the section $s$ can be proven, using a limit argument, 
by observing the family of sections $\{s_i\}_{i\ge 1}$, rather than only the section $s$. 

The best illustration of this phenomenon is the following crucial
observation, which is extremely important in investigating the Grothendieck anabelian section conjecture, and which is du to Tamagawa.

\proclaim{Lemma 2.5}\ Assume that $k$ is finitely generated over the prime field $\Bbb Q$, or that $k$ is a
$p$-adic local field. Let $s:G_k\to \Pi_X$ be a group-theoretic section of the natural projection
$\pr:\Pi_X \twoheadrightarrow G_k$, and $\{X_i[s]\}_{i\ge 1}$ a system of neighbourhoods of the section $s$ (cf. 2.4).
Then the section $s$ is point-theoretic if and only if $X_i[s](k)\neq \varnothing$, for each $i\ge 1$.
\endproclaim

\demo {Proof}\ See [Tamagawa], Proposition 2.8, (iv).
\qed
\enddemo

The above Lemma 2.5 reduces the proof of the Grothendieck anabelian section conjecture, in the case where
$k$ is finitely generated over the prime field $\Bbb Q$, or that $k$ is a
$p$-adic local field, and $\Sigma=\Primes$, to proving the following implication

$$\{\overline {\Sec} _{\Pi_X}\neq \emptyset\} \Longrightarrow \{X(k)\neq \emptyset\}.$$

Thus, at the heart of the Grothendieck anabelian section conjecture, is the following fundamental problem.

\definition{Problem}  How can one produce, under the assumptions of GASC, or the $p$-adic GASC, a rational point
$x\in X(k)$ starting from a group-theoretic section $s:G_k\to \Pi_X$ of $\Pi_X$?
\enddefinition
 
 No systematic approach has yet been developed so far to attack this problem. 
 
 The only situation where a 
 method, or a technique, is available to solve this problem positively is the method, developed by Koenigsmann, and Pop, in the framework of the birational
 version of the p-adic GASC (cf. [Koenigsmann], and [Pop]), and which resorts to a local-global principle for Brauer groups of fields of transcendence 
 degree $1$ over $p$-adic local fields (see the discussion after Theorem 5.2).

\definition{Remarks 2.6} 

{\bf (i)}\ One of the difficulties in investigating the GASC is that, for the time being, one doesn't 
know how to construct sections of arithmetic fundamental groups, and hence test the validity of the conjecture on concrete examples.

In fact, the GASC itself can be viewed as a "rigidity" statement. 

Namely, the only way one knows (so far) to construct sections
of arithmetic fundamental groups of hyperbolic curves, over finitely generated fields of characteristic $0$, is via decomposition groups 
associated to rational points, and these should be the only sections that exist !

One way, however,  to construct such sections is as follows. 

Let $X$ be a proper,
smooth, geometrically connected, hyperbolic algebraic curve over the field $k$, and $\Sigma \subseteq \Primes$ a non-empty set of prime integers.

Consider the exact sequence
$$1\to \Delta_X\to \Pi_X @>{\pr}>> G_k\to 1,\tag {$1.2$}$$
where $\Pi_X$ is the geometrically pro-$\Sigma$ arithmetic fundamental group of $X$. 

Recall the commutative cartesian diagram:
$$
\CD
1 @>>>  \Delta _X     @>>> \Aut (\Delta _X)     @>>> \Out (\Delta _X)  @>>> 1\\
  @.        @A{\id}AA           @AAA              @A{\rho_{X,\Sigma}}AA \\
1 @>>>  \Delta _X       @>>> \Pi_X       @>{\pr}>> G_k   @>>> 1
\endCD
\tag {$1.4$}
$$

In order to construct a continuous group-theoretic section $s:G_k\to \Pi_X$ of the natural projection
$\pr:\Pi_X\twoheadrightarrow G_k$, it is equivalent to construct a continuous homomorphism ${\tilde \rho_{X,\Sigma}}:
G_k\to \Aut (\Delta _X)$, which lifts the
homomorphism  ${\rho_{X,\Sigma}}:G_k\to \Out (\Delta _X)$ above, i.e. such that the following diagram commutes
$$
\CD
G_k@  >{\tilde \rho_{X,\Sigma}}>> \Aut (\Delta _X) \\
@V{\id}VV         @VVV \\
G_k  @>{\rho_{X,\Sigma}}>>  \Out (\Delta _X) \\
\endCD
$$
as follows directly from the fact that the right square in the diagram (1.4) is cartesian.

{\bf (ii)} In light of the Remark (i), it is possible to construct sections $s:G_k\to \Pi_X$ of $\Pi_X$,
if the image of $G_k$ in $\Out (\Delta_X)$ via the natural homomorphism ${\rho_{X,\Sigma}}:G_k\to \Out (\Delta _X)$ 
is a free pro-$\Sigma$ group.

Under this condition one may hope to construct non geometric sections $s$, i.e. sections which do not arise from rational points. 
This is the method used in [Hoshi] to construct non geometric sections in the case where $k$ is a number field, or a $p$-adic local field, and $\Sigma =\{p\}$.

However, this method is unlikely to produce examples of non-geometric sections in the case where $\Sigma=\Primes$. 
Indeed, in this case the homomorphism
${\rho_{X,\Sigma}}:G_k\to \Out (\Delta _X)$ is injective if $X$ is hyperbolic, and $k$ is a number field, or a $p$-adic local field, as is well-known.
\enddefinition

\subhead {2.7}
\endsubhead
One can formulate a birational version of the Grothendieck anabelian section conjecture as follows (see also [Pop]).

There exists a natural exact sequence of absolute Galois groups
$$1\to \Gal (K_X^{\sep}/K_X.\bar k)\to \Gal (K_X^{\sep}/K_X)\to G_k\to 1.$$

Let $\Sigma \subseteq \Primes$ be a non-empty set of prime integers, and
$$\overline G_X\defeq \Gal (K_X^{\sep}/K_X.\bar k)^{\Sigma}$$ 
the maximal pro-$\Sigma$ quotient of the absolute
Galois group $\Gal (K_X^{\sep}/K_X.\bar k)$. 

Let
$$G_X\defeq \Gal (K_X^{\sep}/K_X)/\Ker (\Gal (K_X^{\sep}/K_X.\bar k)\twoheadrightarrow \Gal (K_X^{\sep}/K_X.\bar k)^{\Sigma})$$
be the maximal geometrically pro-$\Sigma$ Galois group of the function field $K_X$. 

Thus, $G_X$ sits naturally
in the following exact sequence
$$1\to \overline G_X\to G_X\to G_k\to 1.$$

Let $x\in X(k)$ be a rational point. Then $x$ determines a decomposition subgroup
$D_x\subset G_X$, which is only defined up to conjugation by the elements of $\overline G_X$, and which maps onto 
$G_k$ via the natural projection $G_X\twoheadrightarrow G_k$. 

More precisely, $D_x$ sits naturally in the following exact sequence
$$1\to M_X\to D_x\to G_k\to 1,$$
where $M_X\isom \hat \Bbb Z (1)^{\Sigma}$ (cf. 3.1, for the definition of the module of roots of unity $M_X$).

The above sequence is known to be split. The set of all possible splittings, i.e. sections
$G_k\to D_x$ of the above exact sequence, is a torsor under $H^1(G_k,M_X)$. The later can be naturally
identified, via Kummer theory, with the $\Sigma$-adic completion $\hat k^{\times, \Sigma}$ of the multiplicative group $k^{\times}$. 

Each section $G_k\to D_x$ of the natural projection $D_x\twoheadrightarrow G_k$
determines naturally a section $G_k\to G_X$ of the natural projection $G_X\twoheadrightarrow G_k$, whose image is contained in $D_x$.

\definition {The Birational Grothendieck Anabelian Section Conjecture (BGASC)} 
\newline
Assume that $k$ is finitely generated over the
prime field $\Bbb Q$, and $\Sigma =\Primes$. Let $s:G_k\to G_X$ be a group-theoretic section of the natural projection
$G_X \twoheadrightarrow G_k$. Then the image $s(G_k)$ is contained in a decomposition subgroup $D_x\subset G_X$ associated to a unique
rational point $x\in X(k)$. In particular, the existence of the section $s$ implies that $X(k)\neq \varnothing$.
\enddefinition

One can, in a similar way, formulate a $p$-adic version of this conjecture.

\definition {A $p$-adic Version of The Birational Grothendieck Anabelian Section Conjecture ($p$-adic BGASC)} 
Assume that $k$ is a $p$-adic local field, i.e. $k$ is a finite extension of $\Bbb Q_p$, 
and $\Sigma =\Primes$. Let $s:G_k\to G_X$ be a group-theoretic section of the natural projection
$G_X \twoheadrightarrow G_k$. Then the image $s(G_k)$ is contained in a decomposition subgroup $D_x\subset G_X$ associated to a unique
rational point $x\in X(k)$. In particular, the existence of the section $s$ implies that $X(k)\neq \varnothing$.
\enddefinition

\subhead
\S 3. Good Sections of Arithmetic Fundamental Groups
\endsubhead
We use the same notations as in $\S1$, and $\S2$. 

We will introduce the notion of (uniformly) good
sections of arithmetic fundamental groups. Point-theoretic sections of arithmetic fundamental groups (cf. Definition 2.2) 
are (uniformly) good sections.

\subhead {3.1}
\endsubhead
Next, we recall the definition of the arithmetic Chern class associated to a group-theoretic section $s$
of the arithmetic fundamental group $\Pi_X$
(cf. [Sa\"\i di], 1.2, and [Esnault-Wittenberg], for more details).

In what follows all scheme cohomology
groups are \'etale cohomology groups. 

First, let ${\hat \Bbb Z}^{\Sigma}$ be the maximal pro-$\Sigma$ quotient of
$\hat \Bbb Z$, and
$$M_X\defeq \Hom (\Bbb Q/\Bbb Z,(K_X^{\sep})^{\times})\otimes_{\hat \Bbb Z}{\hat \Bbb Z^{\Sigma}}.$$

Note that $M_X$ is a free ${\hat \Bbb Z}^{\Sigma}$-module of rank one, and
has a natural structure of $G_k$-module, which is  isomorphic to the $G_k$-module
${\hat \Bbb Z}(1)^{\Sigma}$, where the ``(1)'' denotes a Tate twist, i.e. $G_k$ acts on
${\hat \Bbb Z}(1)^{\Sigma}$ via the $\Sigma$-part of the cyclotomic character. We will refer to $M_X$ as the
module of roots of unity attached to $X$, relative to the set of primes $\Sigma$. 

Let
$$\eta _X^{\diag} \in H^2(X\times X,M_X)$$
be the \'etale Chern class, which is associated to the diagonal embedding $\iota:X\to X\times_kX$,
or alternatively the first Chern class of the line bundle $\Cal O_{X\times X}(\iota (X))$.

There exists a natural identification
(cf. [Mochizuki], Proposition 1.1)
$$H^2(X\times X,M_X)\isom H^2(\Pi_{X\times X} ,M_X).$$

The Chern class $\eta _X^{\diag}$ corresponds via the above identification to an extension class
$$\eta _X^{\diag}\in  H^2(\Pi_{X\times X},M_X).$$
We shall refer to the extension class $\eta _X^{\diag}$
as the extension class of the diagonal. 

Let $s:G_k\to \Pi_X$ be a group-theoretic section of the natural projection $\Pi_X\twoheadrightarrow G_k$. 
Let
$$1\to M_X\to \Cal D\to \Pi_{X\times X}  \to 1\tag {$3.1$}$$
be a group extension, whose class in  $H^2(\Pi_{X\times X},M_X)$ coincides with the extension class
$\eta _X^{\diag}$ of the diagonal. 

By pulling back the group extension $(3.1)$ by the continuous
injective homomorphism
$$(s,\id): G_k\times_{G_k} \Pi_X    \to \Pi_{X\times X},$$
we obtain a natural commutative diagram:

$$
\CD
1@>>>   M_X    @>>>     \Cal D_s   @>>>  G_k\times _{G_k} \Pi_X        @>>> 1\\
  @.        @V{\id}VV           @VVV              @V{(s,\id)} VV \\
1     @>>> M_X  @>>>   \Cal D    @>>>  \Pi_{X\times X}  @>>> 1
\endCD
$$
where the right square is cartesian.

Further, via the natural identification $G_k\times _{G_k}\Pi_X\isom \Pi_X,$ 
the upper group extension $\Cal D_s$ in the above diagram corresponds to a group extension
(which we denote also $\Cal D_s$)
$$1 @>>>  M_X    @>>>     \Cal D_s   @>>>  \Pi_X\to 1.$$ 

We will refer to the class $[\Cal D_s]$
of the extension $\Cal D_s$ in $H^2(\Pi_X,M_X)$ as the extension class
associated to the section $s$. 

\definition{Definition 3.2 (The ($\Sigma$)-\'Etale Chern Class associated to a Section)}
We define the ($\Sigma$)-\'etale Chern class $c(s)\in  H^2(X,M_X)$
associated to the section $s$ as the element of  $H^2(X,M_X)$ corresponding
to the above extension class $[\Cal D_s]$, which is associated to the section $s$,
via the natural identification
$H^2(\Pi_X,M_X)\isom H^2(X,M_X)$.
\enddefinition

Let $\{X_i[s]\}_{i\ge 1}$ be a system of neighbourhoods 
of the section $s$, and $\{\Pi_X[i,s]\}_{i\ge 1}$ the corresponding open subgroups of $\Pi_X$ (cf. 2.4).
Recall that the section $s$ restricts naturally to a group-theoretic section
$$s_i:G_k\to \Pi_X[i,s]$$
of the natural projection $\Pi_X[i,s]\twoheadrightarrow G_k$,
which fits into the following commutative diagram:
$$
\CD
G_k @>{s_i}>> \Pi_X[i,s] \\
@V{\id}VV     @VVV  \\
G_k @>{s}>> \Pi_X
\endCD
$$
where the right vertical homomorphism is the natural inclusion, for each positive integer $i$ (cf. loc. cit).

One can easily observe the following Lemma (See [Sa\"\i di], Lemma 1.3.1, for more details).

\proclaim {Lemma 3.3} For each positive integer $i$, the image of the Chern class
$c(s_{i+1})\in H^2 (X_{i+1}[s],M_X)$
associated to the section $s_{i+1}:G_k\to \Pi_X[i+1,s]$ in $H^2 (X[i,s],M_X)$, via the
corestriction homomorphism
$\cor : H^2 (X_{i+1}[s],M_X)\to H^2 (X_i[s],M_X)$, coincides with
the Chern class $c(s_i)$ associated to the section $s_i:G_k\to \Pi_X[i,s]$.
\endproclaim

\definition{Definition 3.4  (The pro-($\Sigma$)-\'Etale Chern Class associated to a Section)} Let
$\underset{i\ge 1}\to{\varprojlim} H^2 (X_i[s],M_X)$
be the projective limit of the  $H^2 (X_i[s],M_X)$'s, where the transition homomorphisms are the
corestriction homomorphisms. We define the pro-$(\Sigma)$-\'etale Chern class associated to the section $s$,
relative to the system of neighbourhoods $\{X_i[s]\}_{i\ge 1}$, as the element
$\hat c(s)\defeq (c(s_i))_{i\ge 1}\in \underset{i\ge 1}\to{\varprojlim} H^2 (X_i[s],M_X)$
(cf. Lemma 3.3).
\enddefinition

\subhead {3.5}
\endsubhead
Next, we will introduce the notion of (uniformly) good sections of arithmetic fundamental groups.

For each positive $\Sigma$-integer $n$, meaning that $n$ is an integer which is divisible only by primes in
$\Sigma$, the Kummer exact sequence in \'etale topology
$$1\to \mu_n\to \Bbb G_m @>n>> \Bbb G_m\to 1$$
induces naturally, for each positive integer $i$, an exact sequence of abelian groups
$$0 \to \Pic (X_i[s])/n \Pic (X_i[s]) \to  H^2(X_i[s],\mu_{n}) \to _{n} \Br (X_i[s])
\to 0,\tag {$3.2$} $$
which for positive integers $m$ and $n$, with $n$ divides $m$, fits naturally into a commutative diagram:

$$
\CD
0@>>> \Pic (X_i[s])/m \Pic (X_i[s])      @>>>     H^2(X_i[s],\mu_{m})  @>>> _{m} \Br (X_i[s]) @>>> 0\\
  @.        @VVV           @VVV              @VVV \\
0@>>> \Pic (X_i[s])/n \Pic (X_i[s])      @>>>     H^2(X_i[s],\mu_{n})  @>>> _{n} \Br (X_i[s])    @>>> 0
\endCD
$$
where the lower and upper horizontal sequences are the above exact sequence $(3.2)$,
and the vertical homomorphisms are the natural homomorphisms. 

Here,
$\Pic\defeq H^1(\ ,\Bbb G_m)$ denotes the Picard group, $\Br\defeq H^2(\ ,\Bbb G_m)$ is
the Brauer-Grothendieck cohomological group, and for a positive integer $n$: $_{n} \Br\subseteq \Br$ is
the subgroup of $\Br$ which is annihilated by $n$. 

By taking
projective limits, the above diagram induces naturally, for each positive integer $i$, the following exact sequence
$$0\to \underset{n\ \Sigma-\text {integer}}\to{\varprojlim}  \Pic (X_i[s])/n \Pic (X_i[s]) \to H^2(X_i[s],M_X)\to
\underset{n\ \Sigma-\text {integer}}\to{\varprojlim} _n \Br (X_i[s])\to 0.$$

We will denote by
$$\Pic (X_i[s])^{\wedge,\Sigma}\defeq \underset{n\ \Sigma-\text {integer}}
\to{\varprojlim} \Pic (X_i[s])/n \Pic (X_i[s])$$
the $\Sigma$-adic completion of the Picard group $\Pic (X_i[s])$, and
$$T\Br(X_i[s])^{\Sigma}\defeq \underset{n\ \Sigma-\text{integr}}\to{\varprojlim} _{n} \Br (X_i[s])$$
the $\Sigma$-Tate module of the Brauer group $\Br (X_i[s])$. Thus, we have a natural exact sequence
$$0\to\Pic (X_i[s])^{\wedge,\Sigma}\to H^2(X_i[s],M_X)\to T\Br(X_i[s])^{\Sigma}\to 0. \tag {$3.3$}$$

In what follows we will identify $\Pic (X_i[s])^{\wedge,\Sigma}$ with its image in  $H^2(X_i[s],M_X)$
(cf. the above exact sequence (3.3)),
and refer to it as the Picard part of $H^2(X_i[s],M_X)$.

Let 
$$s:G_k\to \Pi_X[1,s]\defeq \Pi_X$$
be a group-theoretic section of $\Pi_X$ as above. For each positive integer $i$,
let 
$$s_i:G_k\to \Pi_X[i,s]$$
be the induced group-theoretic section of $\Pi_X[i,s]$.

By pulling back cohomology classes via the section $s_i$,
and bearing in mind the natural identifications 
$H^2(\Pi_X[i,s],M_X)\isom H^2(X_i[s],M_X)$
(cf. [Mochizuki], Proposition 1.1), we obtain a natural (restriction) homomorphism
$$s_i^{\star}: H^2(X_i[s],M_X)\to H^2(G_k,M_X).$$

Finally, observe that if $k'$ is a finite extension of $k$, and $X_{k'}\defeq X\times _kk'$,
then we have a natural commutative diagram:
$$
\CD
1 @>>>  \Delta _X     @>>> \Pi_{X_{k'}}\defeq \pi_1(X_{k'} ,\eta)^{(\Sigma)}     @>>> G_{k'} @>>> 1\\
  @.        @V{\id}VV           @VVV              @VVV \\
1 @>>>    \Delta _X      @>>> \Pi_X       @>>> G_k   @>>> 1
\endCD
$$
where the right, and middle, vertical arrows are the natural inclusions, and the far right square is cartesian.

In particular, the section $s_k\defeq s:G_k\to \Pi_X$
induces naturally a group-theoretic section
$s_{k'}:G_{k'}\to  \Pi_{X_{k'}}$ of $\Pi_{X_{k'}}$.

\definition{Definition 3.6 (Good and Uniformly Good Sections of
Arithmetic Fundamental Groups)}
We say that the section $s$ is a good group-theoretic section, relative to the system
of neighbourhoods $\{X_i[s]\}_{i\ge 1}$, if the above homomorphisms
$s_i^{\star}: H^2(X_i[s],M_X)\to H^2(G_k,M_X)$,
for each positive integer $i\ge 1$, annihilate the Picard part $\Pic (X_i[s])^{\wedge,\Sigma}$ of
$H^2(X_i[s],M_X)$. In other words,
the section $s$ is good if $\Pic (X_i[s])^{\wedge,\Sigma}\subseteq \Ker {s_i}^{\star}$,
for each positive integer $i$.

We say that the section $s$ is uniformly good, relative to the system
of neighbourhoods $\{X_i[s]\}_{i\ge 1}$, if the induced section
$s_{k'}:G_{k'}\to \Pi_{X_{k'}}$
is good, relative to the system of neighbourhoods of $s_{k'}$ which is naturally induced by the $\{X_i[s]\}_{i\ge 1}$,
for every finite extension $k'/k$.
\enddefinition

It is easy to see that the above definition is independent of the given system of neighbourhoods $\{X_i[s]\}_{i\ge 1}$
of the section $s$. We will refer to a section satisfying the conditions in Definition 3.6 as good, or uniformly good, without necessarily 
specifying a system of neighbourhoods of the section.

The notion of uniformly good sections is motivated by the fact that a necessary condition for a group-theoretic
section $s:G_k\to \Pi_X$ to be point-theoretic, i.e. arises from a $k$-rational point
$x\in X(k)$ (cf. Definition 2.2), is that the section $s$ is uniformly good in the sense of
Definition 3.6 (cf. [Sa\"\i di], Proposition 1.5.2).

\definition{Remarks 3.7}

{\bf (i)} If $k$ is a $p$-adic local field, the conditions of goodness and uniform goodness for the section $s$ are equivalent.

Moreover, if $p\in \Sigma$, the section $s$ is good in this case if and only if $X(k^{\tame})\neq \emptyset$, where $k^{\tame}$ 
is the maximal tame extension of $k$ (cf. [Sa\"\i di], Proposition 1.6.6, and Proposition 1.6.8).

{\bf (ii)}
If $k$ is a number field, one has a local-global principle for (uniform) goodness. Namely, the section $s$ is good in this case
if and only if the section $s_v$ is good, for each place $v$ of $k$ (cf. loc. cit., Proposition 1.8.1). 
\enddefinition

\definition {Remark 3.8}
Another necessary condition for the section $s$ to be point-theoretic, is that the image of the
pro-Chern class $\hat c (s)\in \underset{i\ge 1}\to{\varprojlim}\
H^2(X_i[s],M_X)$ in $\underset{i\ge 1}\to {\varprojlim}\ T\Br (X_i[s])^{\Sigma}$,
via the natural homomorphism
$$\underset{i\ge 1}\to {\varprojlim}\ H^2(X_i[s],M_X)\to \underset{i\ge 1}\to {\varprojlim}\ T\Br (X_i[s])^{\Sigma},$$
equals $0$. 

In other words, if the section $s$ is point-theoretic, then  the associated pro-Chern class
$\hat c (s)$ (cf. Definition 3.4) lies in the Picard part $\underset{i\ge 1}\to{\varprojlim}\ \Pic (X_i[s])^{\wedge,\Sigma}$ of
$\underset{i\ge 1}\to{\varprojlim}\ H^2(X_i[s],M_X)$.

Indeed, Assume that the section $s$ is point-theoretic, meaning that
$s\defeq s_x:G_k\to \Pi_X$
arises from a $k$-rational point $x\in X(k)$ (cf. Definition 2.2). Then there exists a compatible system of rational points
$\{x_i\in X_i[s](k)\}_{i\ge 1}$, i.e. $x_{i+1}$ maps to $x_i$ via the natural morphism 
$X_{i+1}[s]\to X_i[s]$.

For every positive integer $i$, let $\Cal O(x_i)\in \Pic (X_i[s])$
be the degree $1$ line bundle associated to $x_i$. 

The Chern class
$c(s_i) \in H^2(X_i[s],M_X)$, which is associated to the section $s_i$ (cf. Definition 3.2),
coincides with the \'etale Chern class
$c(x_i)\in H^2(X_i[s],M_X)$ associated to the line bundle $\Cal O(x_i)$
(cf. [Mochizuki], Proposition 1.6, and Proposition 1.8 (iii)).
Thus, the pro-Chern class $\hat c(s)$ is a Picard element. More precisely, it 
is the ``pro-Picard element'' induced by the $(\Cal O(x_i))_{i\ge 1}$.

If the above condition is satisfied one says that the Chern class of the section $s$ is algebraic, or that the
section $s$ is well-behaved (cf. [Sa\"\i di], Proposition 1.5.1).

One can easily prove that if the section $s$ is well-behaved, then the section $s$ is good in the sense of Definition 3.6 (cf. [Sa\"\i di], Proposition 1.5.2).
\enddefinition

\subhead
\S 4. The Cuspidalisation of Sections of Arithmetic Fundamental Groups 
\endsubhead

In this section we introduce the problem of cuspidalisation of group-theoretic sections of arithmetic
fundamental groups. 

We follow the notations in $\S 1$. 

In particular, $\Sigma \subseteq \Primes$ is
a non-empty set of prime integers, and we have the natural exact sequence
$$1\to \Delta_X\to \Pi_X @>{\pr_X}>> G_k\to 1,\tag {$1.2$}$$
where $\Pi_X$ is the geometrically pro-$\Sigma$ arithmetic fundamental group of $X$.

\subhead {4.1}
\endsubhead
In this sub-section we recall the definition of (geometrically) cuspidally central, and cuspidally abelian, 
arithmetic fundamental groups of affine hyperbolic curves,
and the definition of cupidally abelian absolute Galois groups of function fields of curves.

\subhead {4.1.1}
\endsubhead
Let $U\subseteq X$ be a non-empty open subscheme of $X$. The geometric point $\eta$ of $X$ determines
naturally a geometric point $\eta$ of $U$, and a geometric point $\bar \eta$ of $\overline U\defeq U\times _k \bar k$.
Write 
$$\Delta_U\defeq \pi_1(\overline U,\bar \eta)^{\Sigma}$$
for the maximal pro-$\Sigma$ quotient of the fundamental group $\pi_1(\overline U,\bar \eta)$ of $\overline U$
with base point $\eta$, and 
$$\Pi_U\defeq  \pi_1(U, \eta)/ \Ker  (\pi_1(\overline U,\bar \eta)\twoheadrightarrow
\pi_1(\overline U,\bar \eta)^{\Sigma})$$ 
for the quotient of  the arithmetic fundamental group
$\pi_1(U, \eta)$ by the kernel of the natural surjective homomorphism
$\pi_1(\overline U,\bar \eta)\twoheadrightarrow \pi_1(\overline U,\bar \eta)^{\Sigma}$, which is a normal subgroup
of $\pi _1(U,\eta)$. 

Thus, we have a natural exact sequence
$$1\to \Delta_U\to \Pi_U @>{\pr_{U,\Sigma}}>> G_k\to 1,$$
which fits into the following commutative diagram:
$$
\CD
1 @>>> \Delta_U  @>>> \Pi_U  @>{\pr_{U,\Sigma}}>> G_k @>>> 1  \\
@.   @VVV     @VVV    @V{\id }VV   \\
1 @>>> \Delta_X   @>>> \Pi_X    @>{\pr_{X,\Sigma}}>> G_k  @>>> 1\\
\endCD
$$
where the left, and middle, vertical homomorphisms are surjective, and are naturally induced by the natural surjective homomorphisms
$\pi_1(\overline U,\bar \eta)\twoheadrightarrow \pi_1(\overline X,\bar \eta)$, and $\pi_1(U,\eta)\twoheadrightarrow \pi_1(X,\eta)$.

Let
$$I_U\defeq \Ker (\Pi_U\twoheadrightarrow \Pi_X)=\Ker (\Delta_U\twoheadrightarrow \Delta_X).$$

We shall refer to $I_U$ as the cuspidal subgroup of $\Pi_U$ (cf. [Mochizuki], Definition 1.5). It is the normal
subgroup of $\Pi_U$ generated by the (pro-$\Sigma$) inertia subgroups at the geometric points of
$S\defeq X\setminus U$.
We have the following natural exact sequence
$$1\to I_U\to \Pi_U\to \Pi_X\to 1.\tag {$4.1$}$$

Let $I_U^{\ab}$ be the maximal abelian quotient of $I_U$. By pushing out the exact sequence $(4.1)$ by the natural
surjective homomorphism $I_U\twoheadrightarrow I_U^{\ab}$, we obtain a natural commutative diagram:
$$
\CD
1 @>>> I_U @>>> \Pi_U @>>> \Pi_X @>>> 1  \\
@. @VVV   @VVV   @V{\id}VV  \\
1 @>>> I_U^{\ab} @>>> \Pi_U^{\c-\ab} @>>> \Pi_X @>>> 1
\endCD
$$

We will refer to the quotient $\Pi_U^{\c-\ab}$ of $\Pi_U$ as the maximal cuspidally abelian quotient of $\Pi_U$,
with respect to the natural homomorphism $\Pi_U\twoheadrightarrow \Pi_X$ (cf. [Mochizuki], Definition 1.5).

Similarly, we can define the maximal cuspidally abelian quotient $\Delta_U^{\c-\ab}$ of $\Delta_U$,
with respect to the natural homomorphism $\Delta_U \twoheadrightarrow \Delta_X$, which sits in a natural exact sequence
$$1 \to I_U^{\ab} \to \Delta_U^{\c-\ab} \to \Delta_X \to 1. \tag{$4.2$}$$

Write $I_U^{\cn}$ for the maximal quotient of  $I_U^{\ab}$ on which the action of $\Delta_X$,
which is naturally deduced from the exact sequence $(4.2)$, is trivial. By pushing out the sequence $(4.2)$
by the natural surjective homomorphism $I_U^{\ab}\twoheadrightarrow I_U^{\cn}$, we obtain a natural exact sequence
$$1\to I_U^{\cn}\to \Delta_U^{\c-\cn}\to \Delta _X\to 1.\tag {$4.3$}$$

Define 
$$\Pi_U^{\c-\cn}\defeq \Pi_U^{\c-\ab}/\Ker (I_U^{\ab}\twoheadrightarrow I_U^{\cn}),$$
which sits naturally in the following exact sequence
$$1\to I_U^{\cn}\to \Pi _U^{\c-\cn}\to \Pi  _X\to 1.\tag {$4.4$}$$

We shall refer to the quotient $\Pi_U^{\c-\cn}$ of $\Pi_U$ as the maximal (geometrically) cuspidally central
quotient of $\Pi_U$, with respect to the natural homomorphism $\Pi_U\twoheadrightarrow \Pi_X$ (cf. loc. cit.).

We have a natural commutative diagram of exact sequences:

$$
\CD
1 @>>> I_U @>>> \Pi_U @>>> \Pi_X @>>> 1\\
@. @VVV     @VVV    @V{\id}VV \\
1 @>>> I_U^{\ab} @>>> \Pi_U^{\c-\ab} @>>> \Pi_X @>>> 1\\
@. @VVV     @VVV    @V{\id}VV \\
1@>>>  I_U^{\cn}  @>>>  \Pi _U^{\c-\cn}  @>>> \Pi  _X @>>> 1\\
\endCD
$$

\subhead {4.1.2}
\endsubhead
Similarly, we have a natural exact sequence of absolute Galois groups
$$1\to G_{\bar k.K_X}\to G_{K_X}\to G_k\to 1,$$
where $G_{\bar k.K_X}\defeq \Gal (K_X^{\sep}/\bar k.K_X)$, and $G_{K_X}\defeq \Gal (K_X^{\sep}/K_X)$.

Let 
$$\overline G_X\defeq G_{\bar k.K_X}^{\Sigma}$$ 
be the maximal pro-$\Sigma$ quotient of $G_{\bar k.K_X}$, and 
$$G_X\defeq  G_{K_X}/ \Ker  ( G_{\bar k.K_X}  \twoheadrightarrow G_{\bar k.K_X}^{\Sigma}),$$
which insert into the following commutative diagram of exact sequences:
$$
\CD
1 @>>> \overline G_{X}  @>>>  G_X  @>{\Tilde {\pr}_{X,\Sigma}}>> G_k @>>> 1  \\
@.   @VVV     @VVV    @V{\id }VV   \\
1 @>>> \Delta_X   @>>> \Pi_X    @>{\pr_{X,\Sigma}}>> G_k   @>>>  1\\
\endCD
$$
where the left vertical maps are the natural surjective homomorphisms.

Let
$$I_X\defeq \Ker (G_{X}\twoheadrightarrow \Pi_X)=\Ker (\overline G_{X}\twoheadrightarrow \Delta_X).$$

We will refer to $I_X$ as the cuspidal subgroup of $G_{X}$. It is the normal subgroup of
$G_{X}$  generated by the (pro-$\Sigma$) inertia subgroups at all geometric closed points of $X$.
We have the following natural exact sequence
$$1\to I_X\to G_{X}  \to \Pi_X\to 1.$$

Let $I_X^{\ab}$ be the maximal abelian quotient of $I_X$. By pushing out the above sequence by the natural
surjective homomorphism $I\twoheadrightarrow I^{\ab}$, we obtain a natural exact sequence
$$1\to I^{\ab}\to G_{X}^{\c-\ab}\to \Pi_X\to 1.\tag {$4.5$}$$

We will refer to the quotient $G_{X} ^{\c-\ab}$ as the maximal cuspidally abelian quotient of $ G_{X}$,
with respect to the natural homomorphism $G_{X}\twoheadrightarrow \Pi_X$.
Note that $ G_{X} ^{\c-\ab}$ is naturally identified with the projective limit
$$\underset{U}\to {\varprojlim}\ \Pi_U^{\c-ab},$$ 
where the limit runs over all open subschemes $U$ of $X$.

\subhead {4.2}
\endsubhead
Next, we consider a continuous group-theoretic section $s:G_k\to \Pi_X$
of the natural projection $\pr_X:\Pi_X\twoheadrightarrow G_k$.

\definition{Definition 4.3 (Lifting of Group-Theoretic Sections)}
Let $U\subseteq X$ be a non-empty open subscheme. We say that a continuous group-theoretic
section $s_U:G_k\to \Pi_U$, of the natural projection  $\pr_U\defeq \pr _{U,\Sigma}:\Pi_U\twoheadrightarrow G_k$
(meaning that $\pr_U\circ s_U=\id _{G_k}$), is a lifting
of the section $s:G_k\to \Pi_X$, if $s_U$ fits into a commutative diagram:
$$
\CD
G_k @>{s_U}>> \Pi_U \\
@V{\id}VV   @VVV \\
G_k @>{s}>>  \Pi _X
\endCD
$$
where the right vertical homomorphism is the natural one. 

More generally, we say that a group-theoretic
section $\tilde s :G_k\to G_{X}$ of the natural projection  $\tilde {\pr}\defeq \tilde {\pr}_{X,\Sigma}:G_{X}\twoheadrightarrow G_k$
(meaning that $\tilde {\pr}\circ \tilde s=\id _{G_k}$) is a lifting of
the section $s$, if $\tilde s$ fits into a commutative diagram:
$$
\CD
G_k @>{\tilde s}>> G_{X} \\
@V{\id}VV   @VVV \\
G_k @>{s}>>  \Pi _X
\endCD
$$
\enddefinition

In connection with the Grothendieck anabelian section conjecture, it is natural to consider the following problem.

\definition {The Cuspidalisation Problem for sections of Arithmetic Fundamental Groups}
Given a group-theoretic section $s:G_k\to \Pi_X$ as above, and a  non-empty open subscheme $U\subseteq X$,
is it possible to construct a lifting $s_U:G_k\to \Pi_U$ of $s$? Similarly, given  a group-theoretic section
$s:G_k\to \Pi_X$ as above, is it possible to construct a lifting $\tilde s : G_k\to G_{X}$ of $s$?
\enddefinition

\definition {Remarks 4.4}

{\bf (i)}
One can easily verify that if the section $s$ is point-theoretic, then the section $s$ can be lifted to a section
$s_U:G_k\to \Pi_U$ of the natural projection $\Pi_U\twoheadrightarrow G_k$, for every open subscheme $U\subseteq X$, and can also be lifted to a section
$\tilde s : G_k\to G_{X}$ of the natural projection $G_X\twoheadrightarrow G_k$.

{\bf (ii)} Note that a positive solution to the cuspidalisation problem, in the case where $k$ is finitely generated over $\Bbb Q$ (resp. $p$-adic local field),
and $\Sigma=\Primes$, plus a positive solution to the BGASC (resp. $p$-adic BGASC), gives a positive solution to the GASC (resp. $p$-adic GASC), as
follows easily from Lemma 2.5.
\enddefinition

Our main result concerning the cuspidalisation problem is the following, which shows that good sections of arithmetic fundamental groups 
behave well with respect to this problem.

Before stating our result, we recall some definitions.

\definition {Definition 4.5}
(i)\ We say that the field $k$ is slim if its absolute Galois group $G_k$ is slim in the sense
of [Mochizuki], $\S0$, meaning that every open subgroup of $G_k$ is centre free. Examples of slim fields include number
fields, and $p$-adic local fields (cf. [Mochizuki2], Theorem 1.1.1). 

One defines in a similar way the notion of a slim profinite group $G$, meaning that
every open subgroup of $G$ is centre free.

(ii)\ We say that the field $k$ is $\Sigma$-regular, if for every prime integer $l\in \Sigma$, and every finite extension $k'/k$,
the $l$-part of the cyclotomic character $\chi _l:G_{k'}\to \Bbb Z_l^{\times}$ is not trivial; or equivalently, if for every prime integer $l\in \Sigma$
the image of the $l$-part of the cyclotomic character $\chi _l:G_k\to \Bbb Z_l^{\times}$ is infinite. 

Examples of $\Sigma$-regular fields (for every non-empty set
$\Sigma$ of prime integers) include number fields, $p$-adic local fields, and finite fields. 

The field $k$ is $\Sigma$-regular if and only if, for every finite extension $k'/k$, the $G_{k'}$-module $M_X$
has no non trivial fixed elements.
\enddefinition

The following is our main result on the cuspidalisation problem (cf. [Sa\"\i di], Theorem 2.6).

\proclaim {Theorem 4.6 (Lifting of Uniformly Good Sections to Cuspidally abelian Arithmetic Fundamental Groups over Slim Fields)} 
Assume that the section $s:G_k\to \Pi_X$ is uniformly good (in the sense of Definition 3.6), and that the field $k$ is slim, and $\Sigma$-regular (cf. Definition 4.5).
Let $U\subseteq X$ be a non-empty open subscheme of $X$, and
$\Pi_{U}^{\c-\ab}$ the maximal cuspidally abelian quotient of $\Pi_{U}$, with respect to the natural surjective
homomorphism $\Pi_{U}\twoheadrightarrow \Pi_X$. 

Then there exists a section $s_U^{\c-\ab}:G_k\to \Pi_{U}^{\c-\ab}$
of the natural projection  $\Pi_{U}^{\c-\ab}\twoheadrightarrow G_k$ which lifts the section $s$, i.e. which inserts into
the following commutative diagram:
$$
\CD
G_k @>s_U^{\c-\ab}>>  \Pi_{U}^{\c-\ab} \\
@V{\id}VV     @VVV  \\
G_k   @>{s}>> \Pi_X
\endCD
$$

Moreover, one can construct for every non-empty open subscheme $U\defeq X\setminus S$ of $X$ a section
$s_U^{\c-\ab}:G_k\to \Pi_{U}^{\c-\ab}$ as above (i.e. which lifts the section $s$), such that for every non-empty open subscheme
$V\defeq X\setminus T$ of $X$, with $U\subseteq V$, we have the following commutative diagram:
$$
\CD
G_k @>s_U^{\c-\cn}>>  \Pi_{U}^{\c-\ab} \\
@V{\id}VV     @VVV  \\
G_k   @>s_V^{\c-\cn}>> \Pi_{V}^{\c-\ab}\\
\endCD
$$
where the right vertical homomorphism is the natural one.
\endproclaim

As a corollary of the above result one obtains the following.

\proclaim {Corollary 4.7 (Lifting of Uniformly Good Sections to Cuspidally abelian Galois Groups over Slim Fields)} Assume that
the field $k$ is slim, and $\Sigma$-regular (cf. Definition 4.5). Let $s:G_k\to \Pi_X$ be a uniformly good group-theoretic section of the natural projection
$\Pi_X\twoheadrightarrow G_k$ (in the sense of Definition 3.6). Then there exists a section $s^{\c-\ab}:G_k\to G_{X}^{\c-\ab}$ of the
natural projection  $G_{X}^{\c-\ab}  \twoheadrightarrow G_k$, which lifts the section $s$,
i.e. which inserts into the following commutative diagram:
$$
\CD
G_k @>s^{\c-\ab}>>   G_{X}^{\c-\ab}\\
@V{\id}VV     @VVV  \\
G_k   @>{s}>> \Pi_X
\endCD
$$
\endproclaim

\subhead {4.8}
\endsubhead
Next, we would like to explain the basic idea behind the proof of Theorem 4.6. See loc. cit. for more details.

Assume that the field $k$ is slim, $\Sigma$-regular, and
the section $s:G_k\to \Pi_X$ is uniformly good (in the sense of Definition 3.6).

Let $U\defeq X\setminus S$ be a non-empty open subscheme of $X$, and
$\Pi_{U}^{\c-\cn}$ the maximal (geometrically) cuspidally central quotient of $\Pi_{U}$, with respect to the natural homomorphism
$\Pi_{U}\twoheadrightarrow \Pi_X$. 

One would like to show, in a first step, that there exists a section $s_U^{\c-\cn}:G_k\to \Pi_{U}^{\c-\cn}$ of the
natural projection  $\Pi_{U}^{\c-\cn}\twoheadrightarrow G_k$ which lifts the section $s$, i.e. which inserts into the following
commutative diagram:
$$
\CD
G_k @>s_U^{\c-\cn}>>  \Pi_{U}^{\c-\cn} \\
@V{\id}VV     @VVV  \\
G_k   @>{s}>> \Pi_X
\endCD
$$

First, one treats the case where the set $S=\{x_i\}_{i=1}^n\subseteq X(k)$ consists of finitely many $k$-rational points.
We will assume, without loss of generality, that $S=\{x\}$ consists of a single rational point $x\in X(k)$.

The maximal (geometrically) cuspidally central quotient
$\Pi_{U_{x}}^{\c-\cn}$ of $\Pi_{U_{x}}$,
with respect to the natural projection
$\Pi_{U_{x}}\twoheadrightarrow \Pi_X$, sits naturally in the following exact sequence
$$1\to M_X\to \Pi_{U_{x}}^{\c-\cn} \to \Pi_X\to 1\tag {$4.6$}$$
(cf. [Mochizuki], Proposition 1.8). 

 By pulling back the group extension $(4.6)$ by the section $s:G_k\to \Pi_X$, we
obtain a group extension
$$1\to M_X\to s^{\star}(\Pi_{U_{x}}^{\c-\cn}) \to G_k\to 1,\tag {$4.7$}$$
which inserts naturally in the following commutative diagram:
$$
\CD
1 @>>> M_X  @>>>  s^{\star}(\Pi_{U_{x}}^{\c-\cn}) @>>> G_k  @>>> 1 \\
@.    @V{\id}VV  @VVV    @V{s}VV   \\
1 @>>> M_X  @>>>  \Pi_{U_{x}}^{\c-\cn} @>>> \Pi_X  @>>> 1
\endCD
$$
where the right square is cartesian. 

The class in $H^2(\Pi_X,M_X)$ of the group
extension $(4.6)$ coincides, via the natural identification  $H^2(\Pi_X,M_X)\isom H^2(X,M_X)$,
with the \'etale Chern class $c(x)\in H^2(\Pi_X,M_X)$ associated to the degree $1$ line bundle $\Cal O(x)$ (cf. [Mochizuki3], lemma 4.2).
The class in $H^2(G_k,M_X)$ of the group
extension $(4.7)$ coincides then with the image $s^{\star} (c(x))$ of the Chern class $c(x)$ via the
(restriction) homomorphism
$s^{\star}:H^2(X,M_X)\to H^2(G_k,M_X)$, which is naturally induced by $s$. 

This image equals $0$, since the
section $s$ is assumed to be good. This follows from the very definition of goodness (cf. Definition 3.6).
Thus, the group extension $(4.7)$ admits group-theoretic splittings. 
A splitting of the group extension $(4.7)$ determines a section $s_U^{\c-\cn}:G_k\to \Pi_{U}^{\c-\cn}$, which lifts the section $s$.
 
The case where the set $S$ consists of finitely many, not necessarily rational, points is treated in a similar way by
using a descent argument, which resorts to the slimness of $k$, and the fact that $k$ is $\Sigma$-regular.

In general, and in order to lift the section $s$ to a section $s_U^{\c-\ab}:G_k\to \Pi_{U}^{\c-\ab}$, one uses the following description of 
$\Pi_{U}^{\c-\ab}$.

For a finite \'etale Galois cover $X'\to X$, with Galois group $\Gal (X'/X)$, let $U'\defeq U\times _XX'$,
and $\Pi_{U'}^{\c-\cn}$ the maximal (geometrically) cuspidally central quotient of $\Pi_{U'}$,
with respect to the natural surjective homomorphism $\Pi_{U'}\twoheadrightarrow \Pi_{X'}$,
which is slim. 

Denote by $\Pi_{U'}^{\c-\cn} \rtimes ^{\out}\Gal (X'/X)$ the profinite group
which is obtained by pulling back the exact sequence 
$$1\to \Pi_{U'}^{\c-\cn}\to \Aut(\Pi_{U'}^{\c-\cn})\to \Out (\Pi_{U'}^{\c-\cn})\to 1,$$ 
by the natural homomorphism:
$$\Gal (X'/X)\to \Out (\Pi_{U'}^{\c-\cn}).$$ 

Thus, we have a natural exact sequence:
$$1\to \Pi_{U'}^{\c-\cn} \to \Pi_{U'}^{\c-\cn} \rtimes ^{\out}\Gal (X'/X)\to \Gal (X'/X)\to 1.$$
which inserts into the following commutative diagram:
$$
\CD
@.   1   @.   1 \\
@.   @VVV  @VVV \\
@.    I_{U'}^{\c-\cn}   @>{\id}>>  I_{U'}^{\c-\cn}   \\
@.  @VVV    @VVV \\
1@>>> \Pi_{U'}^{\c-\cn} @>>> \Pi_{U'}^{\c-\cn} \rtimes ^{\out}\Gal (X'/X) @>>>  \Gal (X'/X) @>>> 1\\
@.   @VVV            @VVV                @V{\id}VV  \\
1@>>> \Pi_{X'}@>>> \Pi_{X} @>>> \Gal (X'/X)@>>> 1\\
@.     @VVV  @VVV   @VVV \\
@. 1 @.   1 @. 1\\
\endCD
$$

Then we have a natural isomorphism:
$$\Pi_{U}^{\c-\ab}\isom \underset{X'\to X}\to {\varprojlim} \Pi_{U'}^{\c-\cn} \rtimes ^{\out}\Gal (X'/X),$$
where the projective limit is taken over all  finite \'etale Galois cover $X'\to X$ (cf. [Sa\"\i di], Proposition 2.5).

\subhead
\S 5. Applications to the Grothendieck Anabelian Section Conjecture
\endsubhead

In this section we state our main applications of the results concerning the cuspidalisation problem for sections of arithmetic fundamental groups, in $\S4$,
to the Grothendieck anabelian section conjecture.

As we already mentioned (cf. Remarks 4.4, (ii)), 
a positive answer to the cuspidalisation problem, in the case where $k$ is finitely generated over $\Bbb Q$ (resp. $p$-adic local field),
and $\Sigma=\Primes$, plus a positive answer to the BGASC (resp. $p$-adic BGASC), implies a positive answer to the GASC (resp. $p$-adic GASC).

The following result of Koenigsmann concerning the $p$-adic BGASC is fundamental (cf. [Koenigsmann]).

\proclaim {Theorem 5.1\ (Koenigsmann)} The $p$-adic BGASC holds true. More precisely, 
assume that $k$ is a $p$-adic local field, and $\Sigma=\Primes$.
Let $s:G_k\to G_X$ be a group-theoretic section of the natural projection
$G_X \twoheadrightarrow G_k$. Then the image $s(G_k)$ is contained in a decomposition subgroup $D_x$ associated to a
unique rational point $x\in X(k)$. In particular, $X(k)\neq \varnothing$.
\endproclaim

This result has been strengthened by Pop, who proved the following (see [Pop]). For a profinite group
$H$, and a prime integer $p$,  we denote by $H''$ the maximal $\Bbb Z/p\Bbb Z$-metabelian quotient of $H$.
Thus, $H''$ is the second quotient of the $\Bbb Z/p\Bbb Z$-derived series of $H$.

\proclaim {Theorem 5.2.\ (Pop)} Assume that $k$ is a $p$-adic local field, which contains a primitive $p$-th root of $1$, and assume $p\in \Sigma$.
Let $s:G_k''\to G_{K_X}''$ be a group-theoretic section of the natural projection
$G_{K_X}'' \twoheadrightarrow G_k''$. Then the image $s(G_k'')$ is contained in a decomposition subgroup $D_x\subset G_{K_X}''$
associated to a unique rational point $x\in X(k)$. In particular, $X(k)\neq \varnothing$.
Here the $(\ \ )''$ of the various profinite groupe
are with respect to the prime $p$, i.e. the second quotients of the $\Bbb Z/p\Bbb Z$-derived series.
\endproclaim

The above Theorem of Pop can be viewed as a very ``minimalistic'' version of the birational Grothendieck anabelian section
conjecture over $p$-adic local fields. Note that the quotient $G_k''$ of $G_k$ is finite in this case. 

Let us mention few words on the proof of the above result of Pop. 

The technical tool in the proof, which produces a rational point $x\in X(k)$ 
starting from a section $s:G_k''\to G_{K_X}''$, is the following. 

For a profinite group $H$, denote by $H'$
the maximal quotient of $H$ which is abelian and annihilated by $p$. Thus, $H'$ is the first quotient of the
$\Bbb Z/p\Bbb Z$-derived series of $H$.  
 
 The existence of the section $s:G_k''\to G_{K_X}''$ implies the existence of a section $s':G_k'\to G_{K_X}'$ of
 the natural projection $G_{K_X}'\twoheadrightarrow G_k'$.  
 
 Let $\Tilde L/K_X$ be the subextension of $K_X^{\sep}/K_X$ with Galois group  $G_{K_X}'$, and $L/K_X$ the subextension of $\Tilde L/K_X$
 corresponding to the subgroup $s'(G_k')$ of $G_{K_X}'$. Then $L$ is a field of transcendence degree $1$ over $k$. For such a field, Pop proved a local-global
 principle for Brauer groups, which generalises a similar principle for function fields of curves over $p$-adic local fields du to Lichtenbaum, and which reads as follows.

The natural diagonal homomorphism
 $$\Br(L)\to \prod _v \Br (L_v),$$
 where the product runs over all rank $1$ valuations $v$ of $L$, is injective. 
 
 Now the existence of the section  $s:G_k''\to G_{K_X}''$ implies that the natural homomorphism
 $$\Br(k)\to \Br (L)$$
 is injective (see [Pop] for more details). 
 
 Let $\alpha\in \Br(k)$ be an element of order $p$. Then $\alpha$ survives in $\Br (L_v)$, i.e. its image is non-zero,  for some rank $1$ valuation
 $v$ of $L$, by the above local-global principle. 
 Pop then proves (the proof is rather technical) that $v$ is the valuation associated to a unique $k$-rational point $x\in X(k)$.

 \subhead {5.3}
 \endsubhead
Using the technique of cuspidalisation of sections of arithmetic fundamental groups, introduced in $\S4$, one can hope to prove the $p$-adic GASC
by reducing it to the $p$-adic version of the BGASC which was proved by Pop. We are able to prove that this is indeed the case, under some additional assumptions.

Before stating our result we will define the notion of good sections of cuspidally abelian absolute Galois groups.

Assume that $k$ is a $p$-adic local field, and $p\in \Sigma$. Let $G_X$ be the geometrically pro-$\Sigma$ absolute Galois group of $K_X$, and 
$G_{X}^{\c-\ab}$ the maximal cuspidally abelian quotient of $G_X$, with respect to the natural surjective homomorphism 
$G_X\twoheadrightarrow \Pi_X$ (cf. 4.1.2). 

Let 
$$\tilde s:G_k\to G_{X}^{\c-\ab}$$ 
be a continuous group-theoretic section of the natural projection
$G_{X}^{\c-\ab}\twoheadrightarrow G_k$. Let $\Tilde L/K_X$ be the subextension of $K_X^{\sep}/K_X$ with Galois group $G_{X}^{\c-\ab}$, 
and $L/K_X$ the subextension of $\Tilde L/K_X$ corresponding to the closed subgroup $\Tilde s(G_k)$ of $G_{X}^{\c-\ab}$.  

We say that the section $\tilde s$ is a good, or tame point-theoretic, group-theoretic section, if the natural homomorphism
$$\Br(k)\to \Br(L)$$ 
is injective (cf. [Sa\"\i di], 1.7, for more details).

Our main result concerning the Grothendieck anabelian section conjecture over $p$-adic local fields is the following.

\proclaim {Theorem 5.4} Assume that $k$ is a $p$-adic local field, and $p\in \Sigma$. Let $s:G_k\to \Pi_X$ be a
group-theoretic section of the natural projection  $\Pi_X \twoheadrightarrow G_k$. Assume that $s$ is a good section
(in the sense of Definition 3.6). Then there exists a section $s^{\c-\ab}:G_k\to G_{X}^{\c-\ab}$ of the
natural projection $G_{X}^{\c-\ab}  \twoheadrightarrow G_k$, which lifts the section $s$.
Furthermore, if the section  $s^{\c-\ab}:G_k\to G_{X}^{\c-\ab}$ is a good section (in the sense of 5.3), 
then  $X(k)\neq \varnothing$.
\endproclaim

\demo{Proof}
The first assertion is Corollary 4.7. 

For a profinite group $H$, and a prime integer $p$, denote by $H'$
the maximal quotient of $H$ which is abelian, and annihilated by $p$. Thus, $H'$ is the first quotient of the
$\Bbb Z/p\Bbb Z$-derived series of $H$. 

The existence of the section  $s^{\c-\ab}:G_k\to G_{X}^{\c-\ab}$
implies the existence of a section $s':G_k'\to  G_X'$ of the natural projection $G_X'\twoheadrightarrow G_k'$. 

Let $\Tilde L/K_X$
be the sub-extension of $K_X^{\sep}/K_X$ which corresponds to the closed subgroup  $\Ker (G_X\twoheadrightarrow G_X')$ of $G_X$, and $L/K_X$
the sub-extension of $\Tilde L/K_X$ which corresponds to the closed subgroup  $s'(G_k')$ of $G_X'$. 

Assume that the section  $s^{\c-\ab}:G_k\to G_{X}^{\c-\ab}$ is good (in the sense of 5.3). Then the natural homomorphism $\Br k\to \Br L$ is injective. 
Under this assumption (which is implied, in the framework of the proof by Pop of Theorem 5.2, by the lifting property
of the section  $s':G_k'\to  G_X'$ to a section $s'':G_k''\to  G_X''$ which is imposed in [Pop]) Pop proves that the image
$s'(G_k')$ is contained in a decomposition subgroup $D_x\subset G_{K_X}'$ associated to a unique rational point $x\in X(k)$ (cf. loc. cit).
In particular,  $X(k)\neq \varnothing$ in this case.
\qed
\enddemo

With the same notations as in theorem 5.4, the author expects that if $s$ is a good section, then any lifting 
 $s^{\c-\ab}:G_k\to G_{X}^{\c-\ab}$ of $s$ (which exists by Corollary 4.7) is automatically good in the sense of $5.3$.
 In particular, following Theorem 5.4, every good section $s$ should be point-theoretic, if $p\in \Sigma$. 
 The author is unable to prove this at the moment of writing 
 this paper. 
 
\subhead {5.5} 
\endsubhead
One can prove, using cuspidalisation techniques of sections of arithmetic fundamental groups, the following (unconditional) version of
the $p$-adic GASC.

\proclaim {Theorem 5.6} Assume that $k$ is a $p$-adic local field, and $\Sigma =\{p\}$. Let $S\subset X$ be a set of closed points of $X$ which is uniformly dense
in $X$ for the $p$-adic topology, meaning that for each finite extension $k'/k$, $S(k')$ is dense in $X(k')$ for the $p$-adic topology. 
Write $\Pi_{X\setminus S}$ for
the geometrically pro-$\Sigma$ quotient of the arithmetic fundamental group $\pi_1(X\setminus S,\eta)$ 
(which is defined in a similar way as $\Pi_X$, cf. 1.1). 

Let $s:G_k\to \Pi_{X\setminus S}$ be a continuous 
group-theoretic section of the natural projection $\Pi_{X\setminus S}\twoheadrightarrow G_k$. Then $s$ is point-theoretic, i.e. the image 
$s(G_k)$ in $\Pi_{X\setminus S}$ is a decomposition group $D_x\subset \Pi_{X\setminus S}$, associated to unique rational point $x\in X(k)$.
\endproclaim

An example of a set $S$ satisfying the assumptions of Theorem 5.6 is the set of algebraic points, in the case where $X$ is defined over a number field

A proof of Theorem 5.6, was communicated orally to the author by A. Tamagawa. The proof relies on the idea of cuspidalisation, and consists in showing that
the section $s$ can be lifted to a section $\Tilde s:G_k\to G_{X}$ of the
natural projection $G_{X}  \twoheadrightarrow G_k$, where $G_X$ is the geometrically pro-$\Sigma$ absolute Galois group of $K_X$. One then reduces the proof to
the $p$-adic version of the BGASC that was proven by Pop.

The proof uses the 
fact that the Galois group $G_k$ of $k$ is topologically finitely generated, and consists in showing that given  a (not necessarily geometrically connected)
ramified Galois cover $Y\to X$ with Galois group $G$, one can "approximate" it by a Galois cover $Y'\to X$ with Galois group $G$, which is ramified 
only above points which are contained in $S$. The later relies on an approximation argument, \`a la Artin, on Hurwitz spaces of covers. 

\subhead {5.7}
\endsubhead
In [Esnault-Wittenberg1] sections of geometrically abelian absolute Galois groups 
of function fields of curves over number fields were investigated.

It is shown in loc. cit. that the existence of such sections implies (in fact is equivalent to) the existence of degree $1$ divisors on the curve, 
under a finiteness condition of the Tate-Shafarevich group of the jacobian of the curve.

As an application of our results on the cuspidalisation of sections of arithmetic fundamental groups, we can prove an analogous result
for good sections of arithmetic fundamental groups. 

Our main result concerning the Grothendieck anabelian section conjecture over number fields is the following.

\proclaim {Theorem 5.8} Assume that $k$ is a number field, and $\Sigma=\Primes$. 
Let $s:G_k\to \Pi_X$ be a group-theoretic section of the natural projection  $\Pi_X \twoheadrightarrow G_k$. 
Assume that $s$ is a uniformly good section in the sense of definition 3.6, and that the jacobian variety
of $X$ has a finite Tate-Shafarevich group. Then there exists a divisor of degree $1$ on $X$.
\endproclaim

\demo{Proof} Follows formally from Corollary 4.7, and Theorem 2.1 in [Esnault-Wittenberg1].
\qed
\enddemo

\subhead
\S 6. On a Weak Form of the $p$-adic Grothendieck  Anabelian section Conjecture
\endsubhead

In this section we discuss a weak form of the $p$-adic GASC.
We will use the following notations.

Let $p>0$ be a fixed prime integer. Let 
$k$ be a $p$-adic local field, i.e. $k$ is a finite extension of $\Bbb Q_p$, $\Cal O_k$ its ring of integers, and $F$ its residue field.

Let $X$ be a proper, smooth,
geometrically connected, and hyperbolic curve over $k$. 

For a non-empty set of prime integers $\Sigma \subseteq  \Primes$, write
$\Pi_X$ for the geometrically pro-$\Sigma$ fundamental group of $X$, which sits in the exact sequence
$$1\to \Delta_X\to \Pi_X\to G_k\to 1,$$ 
where $\Delta _X$ is the maximal pro-$\Sigma$ quotient of the fundamental group of $\overline X\defeq X\times _k\bar k$ (cf. 1.1).

Recall the natural map (cf. 2.1)
$$\varphi_X\defeq \varphi_{X,\Sigma}: X(k)\to  \overline {\Sec} _{\Pi_X}.$$ 

\subhead {6.1}
\endsubhead
Assume that the hyperbolic $k$-curve $X$ has good reduction over $O_k$, i.e. $X$ 
extends to a smooth, proper, and
relative curve $\Cal X$ over $O_k$, and $p\notin \Sigma$. Let $\Cal X_s\defeq \Cal X\times _{O_k}F$ be the special fibre of $\Cal X$.
Let $\xi$ be a geometric point of $\Cal X_s$ above the generic point of $\Cal X_s$.
Then $\xi$ determines naturally an algebraic closure $\overline F$
of $F$, and a geometric point $\bar {\xi}$ of $\overline {\Cal X_s} \defeq \Cal X_s\times _F \overline F$.

There exists a natural exact sequence of profinite groups
$$1\to \pi_1(\overline {\Cal X_s},\bar \xi)\to \pi_1(\Cal X_s, \xi) @>{\pr}>> G_F\to 1.$$

Here $\pi_1(\Cal X_s, \xi)$ denotes the arithmetic \'etale fundamental group of $\Cal X_s$ with base
point $\xi$, $\pi_1(\overline {\Cal X_s},\bar \xi)$ the \'etale fundamental group of $\overline {\Cal X_s}\defeq 
\Cal X_s\times _F \overline F$ with base
point $\bar \xi$, and $G_F\defeq \Gal (\overline F/F)$ the absolute Galois group of $F$.

Write 
$$\Delta_{\Cal X_s}\defeq \pi_1(\overline {\Cal X_s},\bar \xi)^{\Sigma}$$
for the maximal pro-$\Sigma$ quotient of $\pi_1(\overline X,\bar \xi)$,
and
$$\Pi_{\Cal X_s}\defeq  \pi_1(\Cal X_s, \xi)/ \Ker  (\pi_1(\overline {\Cal X_s},\bar \xi)\twoheadrightarrow
\pi_1(\overline {\Cal X_s},\bar \xi)^{\Sigma})$$
for the quotient of  $\pi_1(\Cal X_s, \xi)$ by the kernel of the natural surjective homomorphism
$\pi_1(\overline {\Cal X_s},\bar \xi)\twoheadrightarrow \pi_1(\overline {\Cal X},\bar \xi)^{\Sigma}$,
which is a normal subgroup
of $\pi _1(\Cal X,\xi)$). 

Thus, we have an exact sequence of profinite groups
$$1\to \Delta_{\Cal X_s}\to \Pi_{\Cal X_s} @>{\pr}>> G_F\to 1.\tag {$6.1$}$$

Moreover, after a suitable choice of the base points $\xi$, and $\eta$, there exists a natural commutative specialisation
diagram:

$$
\CD
1 @>>>  \Delta _X     @>>> \Pi_X     @>>> G_k @>>> 1\\
  @.        @VVV           @V{\Sp_X}VV              @VVV \\
1 @>>>    \Delta _ {\Cal X_s}     @>>> \Pi_{\Cal X_s}      @>>> G_F   @>>> 1
\endCD
\tag {$6.2$}
$$
where the left vertical homomorphism $\Sp:\Delta _X \to \Delta _ {\Cal X_s}$ is an 
isomorphism (since we assumed $p\notin \Sigma$),
and the right vertical
homomorphism is the natural projection $G_k\twoheadrightarrow G_F$, as follows easily from the specialisation theory for fundamental groups
of Grothendieck (cf. [SGA1]). 

In fact, In the above commutative diagram (6.2), the right square is cartesian, as follows easily from the slimness
of $\Delta _X$, and the well-known criterion for good reduction of curves (cf. [Sa\"\i di], Lemma 4.2.2).

\subhead {6.2}
\endsubhead
Let $k_X$, and $k_Y$, be two $p$-adic
local fields, i.e. both $k_X$, and $k_Y$, are finite extensions of $\Bbb Q_p$

Let $X$ (resp. $Y$) be a proper, smooth,
geometrically connected, and hyperbolic curve over $k_X$ (resp. $k_Y$). 

Let $\Sigma \subseteq \Primes$ be a non-empty set
of prime integers, and $\Pi_X$ (resp. $\Pi_Y$) the geometrically pro-$\Sigma$ arithmetic fundamental group of $X$ (resp. $Y$),
which sits in the exact sequence $1\to \Delta_X\to \Pi_X\to G_{k_X}\to 1$
(resp. $1\to \Delta_Y\to \Pi_Y\to G_{k_Y}\to 1$). 

Let
$$\sigma :\Pi_X\isom \Pi_Y$$
be an isomorphism between profinite groups.

\proclaim{Lemma 6.3} The isomorphism $\sigma$ fits into a commutative diagram:
$$
\CD
\Delta _X     @>>> \Delta _Y  \\
 @VVV          @VVV    \\
\Pi _X @>{\sigma}>> \Pi_Y   \\
@VVV          @VVV   \\
G_{k_X} @>>>  G_{k_Y} \\
\endCD
\tag {$6.3$}
$$
where the horizontal maps are isomorphisms, which are naturally induced by $\sigma$.

In particular,
the isomorphism $\sigma$ induces naturally a bijection
$$\sigma ^{\sec}:\overline {\Sec} _{\Pi_X}\isom \overline {\Sec} _{\Pi_Y}.$$

Moreover, the natural isomorphism $G_{k_X}\isom G_{k_Y}$ which is induced by $\sigma$ preserves the inertia subgroups,
i.e. maps the inertia subgroup of $G_{k_X}$ isomorphically to the inertia subgroup of  $G_{k_Y}$.
\endproclaim

\demo {Proof} (cf. [Mochizuki4], Proposition 1.2.1, and Lemma 1.3.8).
\qed
\enddemo

\definition {Definition 6.4} We say that the isomorphism $\sigma$ is point-theoretic, if the image of
$X(k)$ in $\overline {\Sec} _{\Pi_Y}$, via the map $\sigma ^{\sec} \circ \varphi_{X,\Sigma} : 
X(k)\to \overline {\Sec} _{\Pi_Y}$,
coincides with $\varphi_{Y,\Sigma} (Y(k_Y))$. In other words $\sigma$ is point-theoretic if it induces naturally a
bijection $\varphi_{X,\Sigma} (X(k_X))\isom \varphi_{Y,\Sigma} (Y(k_Y))$.
\enddefinition

It is natural, in the framework of the $p$-adic GASC, to consider the following question.

\definition {Question 6.5 (A Weak Form of the Grothendieck Anabelian Section Conjecture over $p$-adic Local Fields)}
Let $k_X$, and $k_Y$, be two $p$-adic
local fields. Let $X$ (resp. $Y$) be a proper, smooth,
geometrically connected, and hyperbolic curve over $k_X$ (resp. $k_Y$). 

Assume that $\Sigma=\Primes$.
Let $\Pi_X$ (resp. $\Pi_Y$) be the geometrically pro-$\Sigma$ fundamental group of $X$ (resp. $Y$), which sits in the exact
sequence $1\to \Delta_X\to \Pi_X\to G_{k_X}\to 1$ (resp. $1\to \Delta_Y\to \Pi_Y\to G_{k_Y}\to 1$). Let
$$\sigma :\Pi_X\isom \Pi_Y$$
be an isomorphism between profinite groups. Is $\sigma$ point-theoretic, in the sense of Definition 6.4?
\enddefinition

Note that the validity of the $p$-adic GASC implies a positive answer to Question 6.5.

Although the pro-$\Sigma$ version of the Grothendieck anabelian section conjecture may not hold over $p$-adic
local fields, in the case where $p\notin \Sigma$ (cf. [Sa\"\i di], Proposition 4.2.1), one may ask weather the following
weak form of the pro-$\Sigma$ Grothendieck anabelian section conjecture still holds, if $p\notin \Sigma$.

\definition {Question 6.6}
Let $k_X$, and $k_Y$, be two $p$-adic
local fields. Let $X$ (resp. $Y$) be a proper, smooth,
geometrically connected, and hyperbolic curve over $k_X$ (resp. $k_Y$). 

Let $\Sigma\subset \Primes$ be a non-empty set
of prime integers, with $p\notin \Sigma$.
Let $\Pi_X$ (resp. $\Pi_Y$) be the geometrically pro-$\Sigma$ fundamental group of $X$ (resp. $Y$). Let
$$\sigma :\Pi_X\isom \Pi_Y$$
be an isomorphism between profinite groups. Is $\sigma$ point-theoretic?
\enddefinition

\definition {Remarks 6.7}  A positive answer to Question 6.6 will imply an absolute version of the Grothendieck anabelian conjecture for smooth, proper, and
hyperbolic curves over $p$-adic local fields (cf. [Mochizuki5], Corollary 2.9).
\enddefinition

In connection with Question 6.6, we can prove the following.

\proclaim {Proposition 6.8} Let $k_X$, and $k_Y$, be two $p$-adic
local fields. Let $X$ (resp. $Y$) be a proper, smooth,
geometrically connected, and hyperbolic curve over $k_X$ (resp. $k_Y$). 

Let $\Sigma \subset \Primes$ be a non-empty set of
prime integers, with $p\notin \Sigma$.
Let $\Pi_X$ (resp. $\Pi_Y$) be the geometrically pro-$\Sigma$ fundamental group of $X$ (resp. $Y$). Let
$$\sigma :\Pi_X\isom \Pi_Y$$
be an isomorphism between profinite groups. 

Assume that $X$ (or $Y$) has good reduction over $k_X$ (or $k_Y$), i.e.
$X$ (or $Y$) extends to a proper, and smooth, relative curve over the valuation ring $O_{k_{X}}$ of $k_X$ (or,
over the valuation ring $O_{k_{Y}}$ of $k_Y$). Then $\sigma$ is point-theoretic.
\endproclaim

\demo{Proof} First, $X$ has good reduction over $k_X$ if and only if $Y$ has good reduction over $k_Y$,
as follows easily form the well-known criterion for good reduction, and the last assertion in Lemma 6.3.

Assume that $X$ (resp. $Y$) extends to a
proper and smooth relative curve $\Cal X$ over the valuation ring $O_{k_{X}}$ of $k_X$ (resp. $\Cal Y$ over
the valuation ring $O_{k_{Y}}$ of $k_Y$). Let $\Cal X_s$ (resp. $\Cal Y_s$) be the special fiber of $\Cal X$
(resp. $\Cal Y$). 

Recall the commutative diagrams (6.2):
$$
\CD
1 @>>>  \Delta _X     @>>> \Pi_X     @>>> G_{k_X} @>>> 1\\
  @.        @VVV           @V{\Sp_X}VV              @VVV \\
1 @>>>    \Delta _ {\Cal X_s}     @>>> \Pi_{\Cal X_s}      @>>> G_{F_X}   @>>> 1
\endCD
$$

and

$$
\CD
1 @>>>  \Delta _Y     @>>> \Pi_Y     @>>> G_{k_Y} @>>> 1\\
  @.        @VVV           @V{\Sp_Y}VV              @VVV \\
1 @>>>    \Delta _ {\Cal Y_s}     @>>> \Pi_{\Cal Y_s}      @>>> G_{F_Y}   @>>> 1
\endCD
$$
where $F_X$ (resp. $F_Y$) is the residue field of $k_X$ (resp. $k_Y$).

Let $s=s_x:G_{k_X}\to \Pi_X$ be a group-theoretic section of the natural projection
$\Pi_X\twoheadrightarrow G_{k_X}$, which is point-theoretic, i.e. arises from a
rational point $x\in X(k_X)$. 

Then $s$ induces naturally a group-theoretic section
$s':G_{k_Y}\to \Pi_Y$ of the natural projection $\Pi_Y\twoheadrightarrow G_{k_Y}$, such that the following diagram is commutative:

$$
\CD
G_{k_X} @>>>   G_{k_Y} \\
@V{s}VV             @V{s'}VV \\
\Pi_X  @>{\sigma}>>   \Pi_Y  \\
\endCD
$$
where the upper arrow is the natural isomorphism which is induced by $\sigma$. We will show that
$s'$ is point-theoretic. 

The natural map $\Sp_X \circ s_x:G_{k_X}\to \Pi_{\Cal X_s}$ factorises
as $G_{k_X}\twoheadrightarrow G_{F_X} @>{\bar s}>> \Pi_{\Cal X_s}$, where $\bar s:G_{F_X}\to \Pi_{\Cal X_s}$
is a group-theoretic section of the natural projection $\Pi_{\Cal X_s}\twoheadrightarrow G_{F_X}$.

The section $\bar s$ is point-theoretic since $s$ is, and corresponds to a decomposition subgroup
associated to the point $\bar x\in \Cal X_s (F_X)$, which is the specialisation of the rational point $x\in X(k_X)$.

Let
$$...\subseteq \Delta _{X}[i+1]\subseteq \Delta _{X}[i] \subseteq .....\subseteq \Delta _{X}[1]\defeq \Delta_{X}$$
be a family of open characteristic subgroups of $\Delta_X$, with $\bigcap _{i\ge 1} \Delta _{X}[i]=\{1\}$. We also denote by
$\{\Delta _{X}[i]\}_{i\ge 1}$ the corresponding family of open subgroups of  $\Delta _{\Cal X_s}$,
via the specialisation isomorphism  $\Sp_X: \Delta _X \isom \Delta _ {\Cal X_s}$.

For a positive integer $i$, let $\Pi_{X}[s,i]\defeq \Delta _{X}[i].s(G_{k_X})$, and
$\Pi_{\Cal X_s}[\bar s,i]\defeq \Delta _{X}[i].\bar s(G_{F_X})$. The system of open subgroups
$\{\Pi_{X}[s,i]\}_{i\ge 1}$ (resp. $\{\Pi_{\Cal X_s}[\bar s,i]\}_{i\ge 1}$) correspond to a tower of finite \'etale covers
$$...\to X_{i+1}\to X_i\to ...\to X_1\defeq X,$$
(resp.
$$...\to \Cal X_{s,i+1}\to \Cal X_{s,i}\to ...\to \Cal X_{1,s}\defeq \Cal X_s),$$
where $\{X_i\}_{i\ge 1}$ (resp. $\{\Cal X_{s,i}\}_{i\ge 1}$) form a system of neighbourhoods of the section $s$
(resp. $\bar s$). 

Let $\Pi_{Y}[i]\defeq \sigma (\Pi_{X}[s,i])$, and $\Pi_{\Cal Y_s}[i]\defeq \sigma (\Pi_{\Cal X_s}[\bar s,i])$.
The natural map $\Sp_Y \circ s':G_{k_Y}\to \Pi_{\Cal Y_s}$ factorises
as $G_{k_Y}\twoheadrightarrow G_{F_Y} @>{\bar s'}>> \Pi_{\Cal Y_s}$, where $\bar s':G_{F_Y}\to \Pi_{\Cal Y_s}$
is a group-theoretic section of the natural projection $\Pi_{\Cal Y_s}\twoheadrightarrow G_{F_y}$
(this factorisation is induced by the similar factorisation of $\Sp_X \circ s$).

The $\{\Pi_{Y}[i]\}_{i\ge 1}$ (resp. $\{\Pi_{\Cal Y_{s}[i]}\}_{i\ge 1}$) form a system of neighbourhoods of the section
$s'$ (resp. $\bar s'$), and correspond to a tower of \'etale covers
$$...\to Y_{i+1}\to Y_i\to ...\to Y_1\defeq Y$$ 
(resp. 
$$...\to \Cal Y_{s,i+1}\to \Cal Y_{s,i}\to ...\to \Cal Y_{s,1}\defeq
\Cal Y_s.)$$

The section $\bar s'$ is naturally induced by the section $\bar s$ via the natural isomorphism $\Pi_{\Cal X_k}
\isom \Pi_{\Cal Y_k}$. In particular, the section $\bar s'$ is point-theoretic, since $\bar s$ is point-theoretic,
as follows from the arguments of Tamagawa for characterising decomposition groups of rational points over finite fields
(see for example  [Sa\"\i di-Tamagawa], $\S1$),  hence $\Cal Y_{s,i}(F_{Y})
\neq \varnothing$. 

Also, by the very definition of the system of neighborhouds $\{Y_i\}_{i\ge 1}$, each $Y_i$ has good
reduction over $\Cal O_{k_Y}$, and extends to a smooth and proper model over $\Cal O_{k_Y}$ whose special fibre is
$\Cal Y_{s,i}$. In particular, $Y_i(k_Y)\neq \varnothing$ by the theorem of liftings of smooth points (cf. [SGA1], expos\'e III, Corollaire 3.3).
This implies that the section $s'$ is point-theoretic by Lemma  2.5.
\qed
\enddemo

In connection to question 6.6, and under the condition that the curve $X$ (or $Y$) has potentially good reduction,
one may hope to use Proposition 6.8, plus a descent argument, to give a positive answer to question 6.6 in this case.
The author has no idea, at the moment of writing this paper, on how to perform such a descent argument.

$$\text{References.}$$
\noindent
[Esnault-Wittenberg] Esnault, H., Wittenberg, O., Remarks on cycle classes of sections of the fundamental
group, Mosc. Math. J. 9 (2009), no. 3, 451-467.

\noindent
[Esnault-Wittenberg1] Esnault, H., Wittenberg, O., On abelian birational sections. Journal of the American Mathematical society, 
Volume 23, Number 3, July 2010, Pages 713-724.

\noindent
[Grothendieck]  Grothendieck, A., Brief an G. Faltings, (German), with an english translation on pp. 285-293.
London Math. Soc. Lecture Note Ser., 242, Geometric Galois actions, 1, 49-58, Cambridge Univ. Press,
Cambridge, 1997.

\noindent
[Hoshi] Hoshi, Y., Existence of nongeometric pro-$p$ Galois sections of hyperbolic curves. To appear in Publications of RIMS.

\noindent
[Koenigsmann] Koenigsmann, J., On the section conjecture in anabelian geometry.  J. Reine Angew. Math.  588 
(2005), 221--235.

\noindent
[Mochizuki]  Mochizuki, S., Absolute anabelian cuspidalizations of proper hyperbolic curves.  J. Math. Kyoto
Univ.  47  (2007),  no. 3, 451--539.

\noindent
[Mochizuki1]  Mochizuki, S., The local pro-$p$ anabelian geometry of curves.  Invent. Math.  138  (1999), 
no. 2, 319--423.

\noindent
[Mochizuki2]  Mochizuki, S., Topics surrounding the anabelian geometry of hyperbolic curves. 
Galois groups and fundamental groups,  119--165, Math. Sci. Res. Inst. Publ., 41, Cambridge Univ. Press,
Cambridge, 2003.

\noindent
[Mochizuki3]  Mochizuki, S., Galois sections in absolute anabelian geometry.  Nagoya Math. J. 179 (2005), 17-45.

\noindent
[Mochizuki4]  Mochizuki, S., The absolute anabelian geometry of hyperbolic curves. 
Galois theory and modular forms, 77--122, Dev. Math., 11, Kluwer Acad. Publ., Boston, MA, 2004. 

\noindent
[Mochizuki5]  Mochizuki, S., Topics in absolute anabelian geometry II: decomposition groups and endomorphisms.
Preprint. Available in the home web page of Shinichi Mochizuki.

\noindent
[Pop] Pop. F., On the birational $p$-adic section Conjecture, Compos. Math. 146 (2010), no. 3, 621-637.

\noindent
[Sa\"\i di] Sa\"\i di, M., Good sections of arithmetic fundamental groups. Manuscript.

\noindent
[Sa\"\i di-Tamagawa] Sa\"\i di, M., Tamagawa, A., A prime-to-$p$ version of the Grothendieck anabelian 
conjecture for hyperbolic curves over finite fields of characteristic $p>0$. Publ. Res. Inst. Math. Sci. 45 (2009), no. 1, 135--186. 

\noindent
[SGA1] Grothendieck, A., Rev\^etements \'etales et groupe fondamental. 
S\'eminaire de g\'eom\'etrie alg\'ebrique du Bois Marie 1960--61.

\noindent
[Tamagawa]  Tamagawa, A., The Grothendieck conjecture for affine curves.  Compositio Math.  109. 
(1997),  no. 2, 135--194.

\bigskip

\noindent
Mohamed Sa\"\i di

\noindent
College of Engineering, Mathematics, and Physical Sciences

\noindent
University of Exeter

\noindent
Harrison Building

\noindent
North Park Road

\noindent
EXETER EX4 4QF 

\noindent
United Kingdom

\noindent
M.Saidi\@exeter.ac.uk

\end
\enddocument